\documentclass[%
 aps,pre,
 superscriptaddress,
 sd,%
 amsmath,amssymb,
 reprint,%
]{revtex4-1}
\pdfoutput=1
\usepackage{graphicx}% Include figure files
\usepackage{dcolumn}% Align table columns on decimal point
\usepackage{bm}% bold math
\usepackage{amsmath}
\usepackage{amsfonts}
\usepackage{amssymb}
\usepackage{cancel}
\usepackage{graphicx}
\graphicspath{{Eps_pdf_figs/}}
\usepackage{float}
\usepackage{epstopdf}
\usepackage{wrapfig}
\usepackage{caption}
\usepackage{subcaption}
\usepackage{bigints}
\usepackage{color}
\usepackage{booktabs}
\usepackage{array}

\renewcommand{\d}{\mathrm{d}}

\newcommand{\change}[1]{{\color{black}#1}}

\newcommand{\R}{{\mathbb{R}}}
\newcommand{\s}{{\mathrm{(s)}}}
\newcommand{\sneg}{{(\mathrm{s},0)}}
\newcommand{\spos}{{(\mathrm{s},1)}}
\newcommand{\lin}{{\mathrm{lin}}}
\newcommand{\re}{{\operatorname{Re}}}

\renewcommand{\u}{{\mathrm{(u)}}}
\newcommand{\uneg}{{(\mathrm{u},0)}}
\newcommand{\upos}{{(\mathrm{u},1)}}
\newcommand{\pr}[1]{{\color{black}#1}}
\newcommand{\prjs}[1]{{\color{black}#1}}
\newcommand{\prnew}[1]{{\color{black}#1}}
\newcommand{\review}[1]{{\color{black}#1}}

\draft % marks overfull lines with a black rule on the right

\begin{document}

%\preprint{}

\title{Probability of noise and rate-induced tipping} %Title of paper



\author{Paul Ritchie}
\email{pdlr201@exeter.ac.uk}
\affiliation{Centre for Systems, Dynamics and Control, College of Engineering, Mathematics and Physical Sciences, Harrison Building, University of Exeter, Exeter, EX4 4QF, United Kingdom}
%\homepage[]{http://emps.exeter.ac.uk/mathematics/staff/pdlr201}
%\thanks{}
%\affiliation{University of Exeter}

\author{Jan Sieber}
\email{J.Sieber@exeter.ac.uk}
%\homepage[]{http://empslocal.ex.ac.uk/people/staff/js543/}
%\thanks{}
\affiliation{Centre for Systems, Dynamics and Control, College of Engineering, Mathematics and Physical Sciences, Harrison Building, University of Exeter, Exeter, EX4 4QF, United Kingdom}
\affiliation{EPSRC Centre for Predictive Modelling in Healthcare, University of Exeter, Exeter, EX4 4QJ, United Kingdom}


\date{\today}

\begin{abstract}
  We propose an approximation for the probability of tipping when the
  speed of parameter change and additive white noise interact to cause
  tipping. Our approximation is valid for small to moderate drift
  speeds and helps to estimate the probability of false positives and
  false negatives in early-warning indicators in the case of rate- and
  noise-induced tipping. We illustrate our approximation on a
  prototypical model for rate-induced tipping with additive noise
  using Monte-Carlo simulations. The formula can be extended to close
  encounters of rate-induced tipping and is otherwise applicable to other
  forms of tipping.

  We also provide an asymptotic formula for the critical ramp speed of
  the parameter in the absence of noise for a general class of systems
  undergoing rate-induced tipping.
\end{abstract}

\keywords{Tipping point, rate-induced, noise-induced}%Use showkeys class option if keyword

\maketitle %\maketitle must follow title, authors, abstract and \pacs

%\begin{quotation}
%
%
%\end{quotation}

\section{Introduction}
%\js{Final checks:
%  \begin{itemize}
%  \item Do consistently ``single-mode approximation'' and 'three-mode
%    approximation''.
%  \item Check for consistent american spelling (iza, ize, color) where
%    required.
%  \end{itemize}
%}

The notion of tipping describes the phenomenon
observed in science, where gradual changes to input levels cause a
sudden (in practice possibly catastrophic) change in the
output. Examples of tipping in science include: Arctic sea ice melting
\citep{kwasniok2013predicting}, degradation of coral reefs
\citep{folke2004regime}, dieback of tropical forest and savanna to a
treeless state \citep{hirota2011global} and financial market crashes
\citep{may2008complex}.

Recently, \citet{ashwin2012tipping} attempted to classify the
underlying mechanisms behind any observed tipping event as an example
of either bifurcation-, noise- or rate-induced tipping. The case of a
slow passage through a bifurcation (often a saddle-node), causing a
loss of stability and therefore an abrupt transition to an alternative
stable state \citep{lenton2013environmental}, is called
\emph{bifurcation-induced tipping}. In contrast, random (rare) jumps
between attractors of an underlying deterministic system due to fluctuations is classified as \emph{noise-induced
  tipping}. \emph{Rate-induced tipping} occurs when a system fails to
track the continuously changing quasi-steady state
\citep{ashwin2012tipping} because the parameter drift speed exceeds a
certain critical rate. For more general definitions and properties of
bifurcation- or rate-induced tipping we refer to
\citet{ashwin2015parameter}.

A research area related to tipping is the study of  generic early-warning indicators \citep{scheffer2009early}. Increase of autocorrelation and variance in output time series are two statistical indicators which are based on the phenomenon commonly known as `critical slowing down' as a system parameter approaches a bifurcation value \citep{van2014critical}. The idea is that far from a bifurcation, the state of the system behaves like  an overdamped particle in a slowly softening potential well \citep{ritchie2015early}. If a small perturbation is made to the particle there will be a fast recovery back to the equilibrium \citep{scheffer2012anticipating}. However, as the bifurcation is approached, the well softens and the recovery from a small perturbation will be slower such that one generically observes an increase in the autocorrelation and variance in output time series \citep{lenton2012early}. On the other hand, for purely noise-induced transitions no bifurcation point is approached and therefore there is debate into the usefulness of the early-warning signals for this type of tipping \citep{dakos2015resilience,PhysRevE.93.032404}. Rate-induced tipping does not involve a loss of stability \citep{perryman2014adapting} and therefore \citet{ashwin2012tipping} commented that there is no reason to suggest the early-warning indicators should be present. However, it has been shown that for a prototypical model for rate-induced tipping the autocorrelation and variance increase before the closest encounter with the critical rate occurs \citep{ritchie2015early}.  

% Examples from climate science, which show evidence of critical slowing down and thus, have been classed as bifurcation-induced tipping events include: desertification of North Africa and the end of the Younger Dryas and greenhouse Earth \citep{dakos2008slowing}. Whereas, the lack of early-warning indicators for the Dansgaard-Oeschger events, imply this is an example of a noise-induced transition \citep{ditlevsen2010tipping}.

In the study of palaeoclimate records \cite{dakos2008slowing,ditlevsen2010tipping} early-warning indicators have been tested on events in the past when tipping has been known to occur. However, testing the early-warning indicators against historical examples is susceptible to statistical mistakes as one selects data conditioned on the system having tipped \citep{boettiger2012early}. 
%\js{The paragraph below is unclear: you say there is little literature but then cite 2 papers. What do they do? Do they agree with your statement that there is little literature? Or do they research this topic. $\Rightarrow$} 

A natural progression will be to use the early-warning indicators to try and predict future tipping events. Though, this \pr{raises} such questions as, if we were to observe an increase in both the autocorrelation and variance of a time series does this mean that the system will tip? \prjs{\citet{boettiger2013no} show there is an increased rate of false positives in early-warning indicators for simulated systems that experience transitions purely by chance. Furthermore, \citet{drake2013early} suggests that stochastic switching can be anticipated but argues that any statistics to be used as early-warning requires decision theory to balance the strength of evidence against the cost and benefits of early-warnings and false positives.}  

% There currently exists little literature that addresses false negatives and false positives created by the early-warning indicators, \prjs{though some attempts have been made} for noise-induced transitions \citep{drake2013early,boettiger2013no}. %\js{$\Leftarrow$}

%\js{I would drop this as we said something similar above? $\Rightarrow$} It has also been eluded to but not explicitly stated that for an example of rate-induced tipping with additive noise the early-warning indicators were shown to be present but the system was not certain to tip \citep{ritchie2015early}. \js{$\Leftarrow$}

This paper provides a generic approximation for the probability of a
prototypical model for rate-induced tipping with additive white
noise. It is structured as follows: \change{Section
  \ref{sec:deterministic} gives a general asymptotic approximation for
  the critical rate for a class of deterministic systems with
  rate-induced tipping.} Section \ref{sec: Methods} derives the
approximation for the probability \change{of a noise-induced escape
  during ramp of a system parameter that does not quite reach the
  critical rate}. In Section \ref{sec: Setup}, we illustrate the
general approximation result with the prototype model introduced by
\citet{ashwin2012tipping}, which we then systematically study in
dependence of its parameter in Section \ref{sec:
  Probability}. Finally, Section \ref{sec: Discussion} the limits of
our approximations and how further developments can address these.

\section{Rate-induced tipping in systems with a ramped parameter}
\label{sec:deterministic}
A general scenario for the phenomenon of rate-induced tipping
was considered by \citet{ashwin2012tipping,ashwin2015parameter}. Assume
that a parameter $\lambda$ corresponds to a shift of the coordinate
system:
\begin{equation}
  \label{eq:det}
  \dot x=f(x+b\lambda)\mbox{,\quad}x(t)\in\mathbb{R}^n\mbox{,\ }
  b\in\mathbb{R}^n\mbox{,\ }\lambda\in\mathbb{R}\mbox{,}
\end{equation}
where the vector $b$ is the direction of the shift and $\lambda$ is
the (scalar) amount. For each fixed $\lambda$ the stability of, for
example, equilibria of \eqref{eq:det} is identical. However, when
$\lambda$ is time-dependent, then there can be \emph{critical rates}
\citep{wieczorek2011excitability} of change of $\lambda$.

\paragraph{Linear shift}
The simplest example discussed in \citet{ashwin2012tipping} is a
linear parameter shift, that is, $\lambda=r_\lin t$ (with $r_\lin>0$
constant). One of the cases studied in \citet{ashwin2012tipping} was
assuming that the system in co-moving coordinates
\begin{equation}\label{eq:linshift}
\dot y=f(y)+r_\lin b \mbox{\quad(where $y=x+br_\lin t$)}
\end{equation}
has a saddle-node bifurcation at $r_\lin=r_0>0$, $y=y_0$ with a stable branch \pr{$y^\s[r_\lin]$} and
an unstable branch \pr{$y^\u[r_\lin]$} of equilibria emerging for
$r_\lin\in[0,r_0)$. \pr{(We will be using square brackets to denote
  branches of equilibria to avoid confusion with time dependence.)}
These equilibria for $y$ correspond to stable and unstable invariant
lines \pr{$x^\s(t)=y^\s[r_\lin]-br_\lin t$,
  $x^\u(t)=y^\u[r_\lin]-br_\lin t$} of the original system
\eqref{eq:det} for $r_\lin\in[0,r_0)$. For $r_\lin<r_0$, all initial
conditions $x(0)$ near $y^\s[r_\lin]$ follow $x^\s(t)$ for all $t>0$,
while for $r_\lin>r_0$ this invariant line no longer exists such that
the rate $r_0$ is critical. This scenario corresponds to a saddle-node
bifurcation in the co-moving coordinates \eqref{eq:linshift} using
$y$. Increasing $r_\lin$ gradually corresponds to a slow passage
through a saddle-node bifurcation.

\paragraph{Ramped shift}  A more complex scenario is the case
where $\lambda$ is ``ramped up'', that is, $\lambda\to0$ for $t\to-\infty$,
$\lambda\to\lambda_{\max}$ for $t\to+\infty$, and $\dot \lambda(t)>0$ for all
$t$. A prototype system for this ramping scenario was studied in
\cite{ashwin2012tipping,ritchie2015early}, and is used in
Section~\ref{sec: Setup} for illustration. For the one-dimensional
case \citet{ashwin2015parameter} gave topological criteria (and a
general definition) for rate-induced tipping with a ramped parameter
$\lambda$. For the general case \eqref{eq:det} we assume that the change
of the ramp in $\lambda$ is itself given by a scalar differential
equation. Define
\begin{align*}
  r&=\max\{\dot \lambda(t):t\in\mathbb{R}\} &&\mbox{(maximal ramp speed)}\\
  \epsilon&=r/\lambda_{\max} &&\mbox{(sharpness of ramp)}\\
  \mu(t)&=\lambda(t)/\lambda_{\max}&&\mbox{($\lambda$ normalized to $[0,1]$)},
\end{align*}
then, assuming $\dot{\lambda}$ is bounded and using the new parameters $r$ and $\epsilon$, $x$ and $\mu$ are the solution of an autonomous extended system:
\begin{align}
    \label{eq:xramp}
    \dot x&=f\left(x+b\frac{r}{\epsilon}\mu\right)\mbox{,}\\
    \dot \mu&=\epsilon\Gamma(\mu)\label{eq:lramp}\mbox{,}
\end{align}
where $\Gamma$ normalizes $\dot{\lambda}$ to $[0,1]$. Since $\dot\lambda$ is always positive, $\Gamma$ satisfies the following properties
\begin{equation}\label{eq:Lambda:ass}
  \begin{split}
    \Gamma(0)&=\Gamma(1)=0\mbox{,\
    }\max\{\Gamma(\mu):\mu\in[0,1]\}=1\mbox{,\ and\ }\\
    \Gamma(\mu)&>0\mbox{\ for all $\mu\in(0,1)$.}
  \end{split}
\end{equation}
Let us also assume that $\lambda$ approaches its limits at an
exponential rate such that $\Gamma'(0)>0$, $\Gamma'(1)<0$, and that
$\Gamma$ is only equal to $1$ in a single point
$\mu_\mathrm{crit}\in(0,1)$ and that $\Gamma''(\mu_\mathrm{crit})<0$.

If system \eqref{eq:linshift} has a saddle-node bifurcation at
$y=y_0$, \pr{$r:=r_\lin=r_0$}, connecting a stable branch \pr{$y^\s[r]$} of
equilibria of \eqref{eq:linshift} and a branch \pr{$y^\u[r]$} with a single
degree of instability for $r\in[0,r_0]$, then we can make the
following statement about the existence of a critical rate
$r_c(\epsilon)$ for sufficiently small $\epsilon$.

The system~\eqref{eq:xramp}--\eqref{eq:lramp} has (at least) $4$
equilibria:
\begin{itemize}
\item \pr{$x_\mathrm{eq}^\sneg:=y^\s[r]\bigg|_{r=0}$}, $\mu=0$ with one unstable direction,
\item \pr{$x_\mathrm{eq}^\uneg:=y^\u[r]\bigg|_{r=0}$}, $\mu=0$ with two
  unstable directions,
\item \pr{$x_\mathrm{eq}^\spos:=y^\s[r]\bigg|_{r=0}-br/\epsilon$}, $\mu=1$ (stable),
\item \pr{$x_\mathrm{eq}^\upos:=y^\u[r]\bigg|_{r=0}-br/\epsilon$}, $\mu=1$ with one unstable
  direction.
\end{itemize}
For sufficiently small $\epsilon$ there are three possible scenarios
for system~\eqref{eq:xramp}--\eqref{eq:lramp} depending on $r$ and a
critical rate $r_c(\epsilon)$, illustrated in Figure~\ref{Phase
  planes}.

\begin{figure}[h!]
        \centering
        \subcaptionbox{\label{Phase plane mux}}[0.45\linewidth]
                {\includegraphics[scale = 0.3]{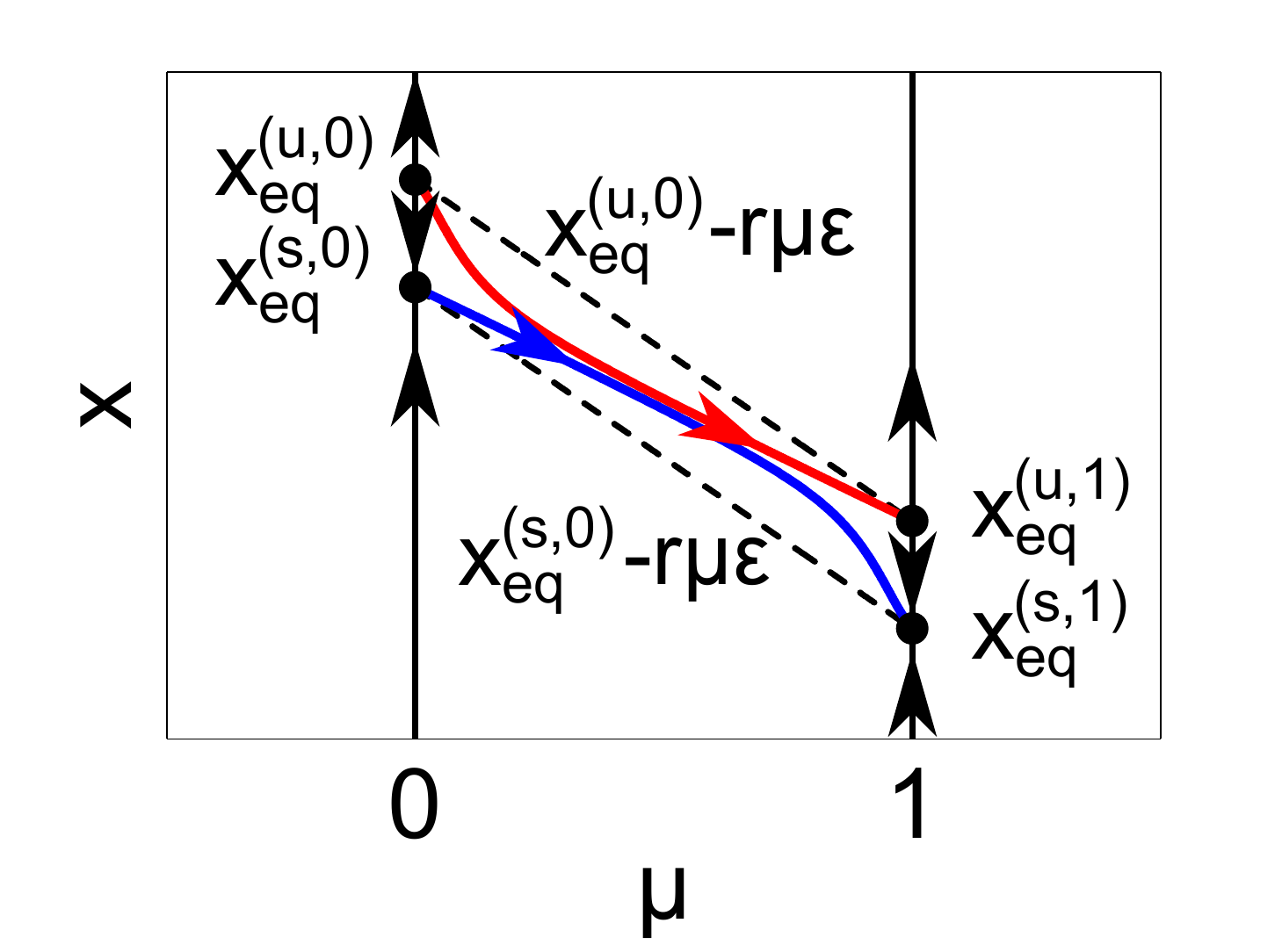}}
        \hfill
        \subcaptionbox{\label{Phase plane muy}}[0.45\linewidth]
                {\includegraphics[scale = 0.3]{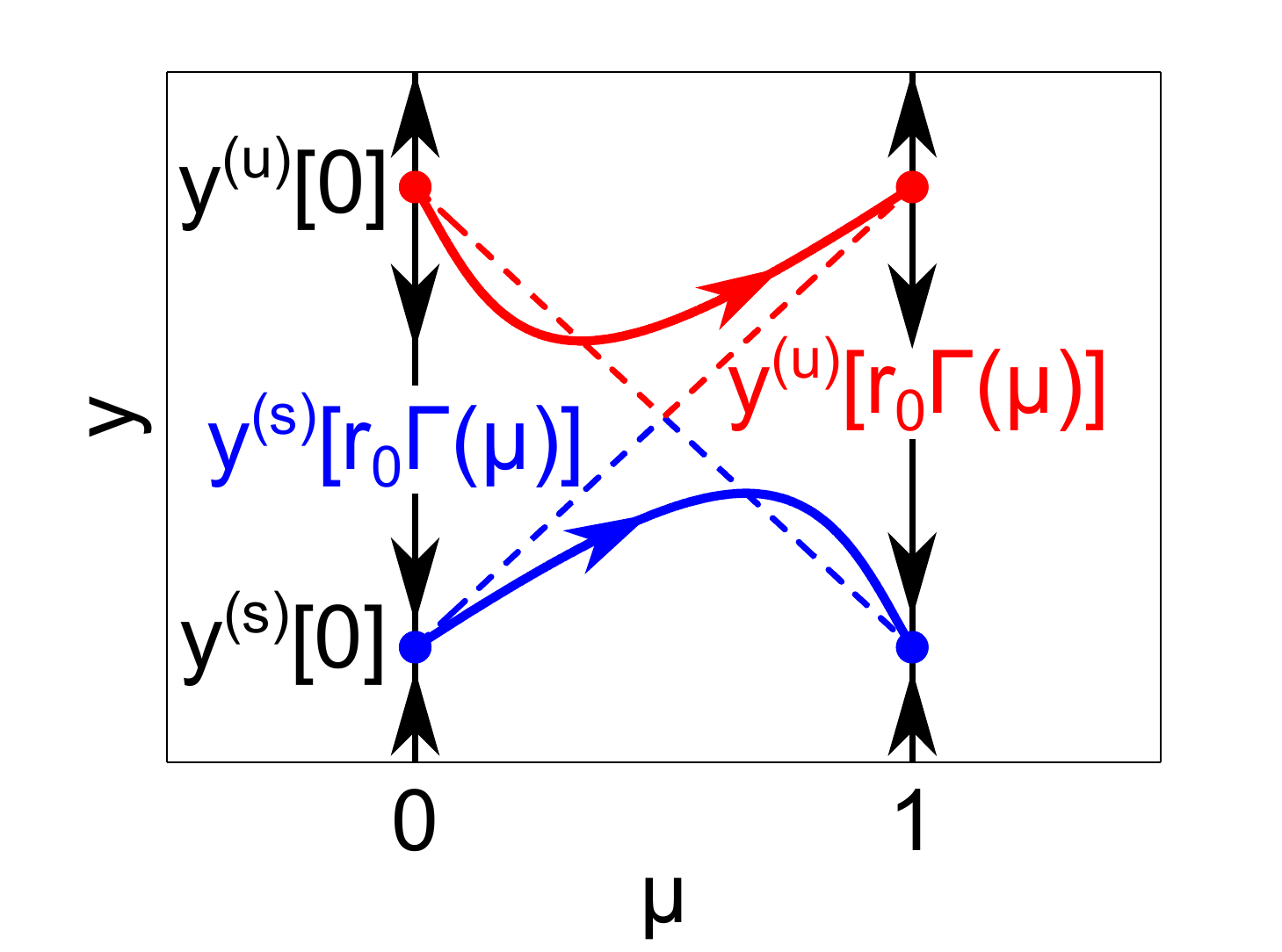}}
        ~ 
        \caption{\pr{(a) Phase plane of system \eqref{eq:xramp}--\eqref{eq:lramp} for the scenario of tracking, $r<r_c(\epsilon)$. Black dashed lines are the stable (lower) and unstable (upper) branches of equilibria in the limit $\epsilon = 0$. Solid blue and red curves represent conecting orbits between $(x,\mu) = (x_\mathrm{eq}^\sneg,0)$ and $(x_\mathrm{eq}^\spos,1)$ and $(x_\mathrm{eq}^\uneg,0)$ to $(x_\mathrm{eq}^\upos,1)$ respectively. (b) Phase plane of system \eqref{scaled y ode}--\eqref{mu ramp} for $r = r_0$, tracking scenario for $\epsilon>0$. The blue (lower) dashed line displays the branch of stable equilibria and the red (upper) dashed line is the branch of unstable equilibria in the limit $\epsilon = 0$. Connecting orbits between $(y,\mu) = (y^\s[0],0)$ and $(y^\s[0],1)$ and $(y^\u[0],0)$ and $(y^\u[0],1)$ given by solid blue and red curves respectively. Parameters: $b=1$, $\epsilon = 0.21$.}}\label{Phase planes}
\end{figure}

\textbf{Tracking:} $r<r_c(\epsilon)$. There is a connecting orbit from
$(x_\mathrm{eq}^\sneg,0)$ to $(x_\mathrm{eq}^\spos,1)$. In this case,
solutions $(x(t),\mu(t))$ starting close to $(x_\mathrm{eq}^\sneg,0)$
stay close to $(x_\mathrm{eq}^\sneg-br\mu(t)/\epsilon,\mu(t))$ for all
$t$ (the distance goes to $0$ as $\epsilon\to0$), see Figure~\ref{Phase plane mux}.

\textbf{Critical:} $r=r_c(\epsilon)=r_0+O(\epsilon)$. There is a
saddle-to-saddle connection from $(x_\mathrm{eq}^\sneg,0)$ to
$(x_\mathrm{eq}^\upos,1)$ in
system~\eqref{eq:xramp}--\eqref{eq:lramp}. The first-order expansion
for $r_c$ in $\epsilon$ is 
\begin{equation}\label{eq:criticalrate}
  r_c(\epsilon)=r_0+\epsilon \sqrt{\frac{-r_0\Gamma''(\mu_{\mathrm{crit}})}{2a_0a_2}} +O(\epsilon^2)\mbox{,}
\end{equation}
where $a_0=w_0^Tb$, $a_2=\frac{1}{2}w_0^T\partial^2f(y_0)v_0^2$, and $w_0$ and
$v_0$ are the left and right nullvectors of $\partial f(y_0)$, scaled
such that $w_0^Tv_0=1$ and $a_0a_2>0$. The coefficients $a_0$ and
$a_2$ are the expansion coefficients when one inserts $y=y_0+v_0z$
and $r_\lin=r\Gamma(\mu)-r_0$ into \eqref{eq:linshift}, applies
$w_0^T$, and truncates to second-order terms:
\begin{equation}
  \dot
  z=a_0(r\Gamma(\mu)-r_0)+a_2z^2+O(z^3)\mbox{.}
  \label{ODE rc expansion}
\end{equation}
See \review{Appendix \ref{app:crit rate expansion}} for
details of the derivation of the first-order expansion for $r_c$ given
by equation \eqref{eq:criticalrate}.

\textbf{Escape:} $r>r_c(\epsilon)$. There are initial conditions
for \eqref{eq:xramp}--\eqref{eq:lramp} arbitrarily close to
$(x_\mathrm{eq}^\sneg,0)$ that escape, following the unstable manifold of
$(x_\mathrm{eq}^\upos,1)$.

The expression for $r_c(\epsilon)$ shows that in the limit $\max\dot
\lambda\ll\lambda_{\max}$ ($\epsilon\ll1$, long ``gentle'' ramps) the
critical rate $r_c(\epsilon)$ for the ramp approaches the critical
rate $r_0$ for the linear shift from above.

The expansion \eqref{eq:criticalrate} for the saddle-connection is
determined entirely by quantities close to the saddle-node, because
$x(t)$ is $\epsilon$-close to \pr{$y^\s[r\Gamma(\mu)]$} for
$\mu<\mu_\mathrm{crit}$, and it is $\epsilon$-close to
\pr{$y^\u[r\Gamma(\mu)]$} for $\mu>\mu_\mathrm{crit}$ (recall
that $\Gamma$ equals $1$ only for $\mu=\mu_\mathrm{crit}$).

% In the limit $\epsilon=0$ there exists a two-dimensional invariant manifold for
% the system in co-moving coordinates ($y=x+br\mu/\epsilon$)
% \begin{align}
% \label{scaled y ode}
% \dot y&=f(y)+br\Gamma(\mu), \\
% \dot\mu&=\epsilon\Gamma(\mu) &&\mbox{(equals $0$ for $\epsilon=0$).}
% \label{mu ramp}
% \end{align}  
% The saddle-to-saddle connection and the critical rate $r_c(\epsilon)$
% exists for those $\epsilon$ for which this invariant manifold persists
% and the unstable manifold of $(x,\mu)=(x_\eq^\sneg,0)$ enters its
% basin of attraction.

Figure~\ref{Phase plane muy} illustrates the phase space for the
  example from Section~\ref{sec: Setup} in co-moving
  coordinates (for scalar  $y$)
  \begin{align}
    \label{scaled y ode}
    \dot y&=f(y)+br\Gamma(\mu), \\
    \dot\mu&=\epsilon\Gamma(\mu)
    \label{mu ramp}
  \end{align}  
  with $r=r_0$. This is an illustration of the tracking scenario
  $(r<r_c(\epsilon))$ and so there is a connecting orbit between
  $(y,\mu)=(y^\s[0],0)$ and $(y^\s[0],1)$. The distance of the
  connecting orbit to $y^\s[r_0\Gamma(\mu)]\to 0$ as $\epsilon\to 0$.

The branches of equilibria $y^\s[r\Gamma(\mu)]$ and $y^\u[r\Gamma(\mu)]$ change their arrangement depending on the value of $r$ in relation to $r_0$. For $r<r_0$ there exist a continuous branch of stable (unstable) equilibria connecting $y^\s[0]$ $(y^\u[0])$ between $\mu = 0$ and $\mu = 1$. For $r=r_0$ the two branches meet as depicted by Figure~\ref{Phase plane muy}. The connections then break up for $r>r_0$ such that two separate saddle-node bifurcations are formed. \review{Appendix \ref{app:phase planes}} contains illustrations for the phase space for other rates.

\section{Noise-induced escape during ramp near but below critical rate}
\label{sec: Methods}
\setcounter{paragraph}{0}
We consider the effect of additive noise for the scalar
setting. Then, system~\eqref{eq:xramp}--\eqref{eq:lramp} changes into
a scalar stochastic differential equation (SDE) for a random variable
$X_{t}$
 \begin{align}
\mathrm{d}X_{t} &= f\left(X_{t}+\frac{r}{\epsilon}\mu\right)\mathrm{d}t + \sqrt{2D}\mathrm{d}W_{t}\mbox{,}
\label{gen SDE}\\
\dot\mu&=\epsilon\Gamma(\mu)\mbox{,}\label{eq:lramp:2}
\end{align}
where $W_{t}$ is standard Brownian motion, the intensity of the noise
is given by $\sqrt{2D}$, and $D$ is a constant diffusion coefficient.
We assume that the deterministic part is as described in
Section~\ref{sec:deterministic}. The deterministic part in \eqref{gen
  SDE} corresponds to a choice of $b$ equal to $1$ in the general
equation~\eqref{eq:xramp}. \prjs{Setting $Y_t=X_t+r\mu/\epsilon$ gives}
%\js{\textbf{There is a problem here:} \eqref{eq:linshift1d} has a
%  $t$-dependent $\mu$, so speaking of a saddle-node doesn't make
%  sense. I think $\dot y=f(y)+r$ has a saddle-node at $r=r_0$.)}
\begin{align}
  \label{eq:linshift1d}
  \d Y_t&=(f(Y_t)+\pr{r}\Gamma(\mu(t)))\d t+\sqrt{2D}\mathrm{d}W_{t} 
\end{align}
where $\dot y=f(y)+r\Gamma(\mu)$ follows the scenario from Section~\ref{sec:deterministic} when treating $\mu$ as a parameter. That is, for $r=r_0$, the system touches a saddle-node non-transversally when $\mu$ crosses $\mu_{\mathrm{crit}}$. As introduced in Section~\ref{sec:deterministic}, the $y^\s$ and $y^\u$ are the stable and unstable branches. They are arranged such that
  $y^\s[r\Gamma(\mu)]<y^\u[r\Gamma(\mu)]$, as shown in
  Figure~\ref{Phase plane muy}, when
  $r\Gamma(\mu)\in[0,r_0)$ (consistent with our setup in
Section~\ref{sec:deterministic} and the example in Section~\ref{sec:
  Setup}). This means that
\begin{displaymath}
  x_\mathrm{eq}^\sneg=y^\s[r\Gamma(\mu)]\bigg|_{r=0}<x_\mathrm{eq}^\uneg=y^\u[r\Gamma(\mu)]\bigg|_{r=0}
\end{displaymath}
are equilibria of the deterministic part of \eqref{gen SDE} combined
with $\dot\mu=\epsilon\Gamma(\mu)$ (identical to \eqref{eq:lramp} with
the same assumptions \eqref{eq:Lambda:ass} on $\Gamma$).  Expressions
below also use the potential $U_r$ for the deterministic
part of \eqref{eq:linshift1d}:
\begin{displaymath}
  U_r(y)=-\smallint f(y)\d y-r\Gamma y\mbox{.}
\end{displaymath}
% \paragraph{Deterministic case \textup{(}$D=0$\textup{):}} for
% $|\epsilon|\ll1$ the trajectory $X_t$ follows the path
% \begin{equation}\label{eq:xu}
%   \tilde  x(t)=x_\mathrm{eq}^\s-\mu(t)-\epsilon\dot\mu(t)/f'(x_\mathrm{eq}^\s)+O(\epsilon^2)\mbox{,}
% \end{equation}
% which is a heteroclinic
% connection between the saddle $x_\mathrm{eq}^\s-\mu_-$ and the sink
% $x_\mathrm{eq}^\s-\mu_+$ of \eqref{gen SDE}--\eqref{eq:genlambda}
% with $D=0$.

% If $\mu_+-\mu_->x_\mathrm{eq}^\u-x_\mathrm{eq}^\s$ then there
% exists a critical ramping speed $\epsilon_c>0$ for which the saddles
% $x_\mathrm{eq}^\s-\mu_-$ and $x_\mathrm{eq}^\u-\mu_+$ of
% \eqref{gen SDE}--\eqref{eq:genlambda} with $D=0$ are connected (see
% \cite{ashwin2015parameter}). If $\epsilon<\epsilon_c$ then
% trajectories starting from $x_\mathrm{eq}^\s-\mu_-$ at
% $\mu\approx \mu_-$ do not escape during the ramp (that is,
% they converge to $x_\mathrm{eq}^\s-\mu_+$). If
% $\epsilon>\epsilon_c$ then trajectories starting from
% $x_\mathrm{eq}^\s-\mu_-$ escape during the ramp beyond the
% equilibrium $x_\mathrm{eq}^\u-\mu_+$ (called \emph{rate-induced
%   tipping} in \cite{ashwin2015parameter}).

\paragraph*{Stationary case with noise } (Equation\,\eqref{gen SDE}
with $r=0$, which is independent of $\mu$) For $D>0$ there will be a fixed
escape rate $\kappa$ from the basin of attraction of
\pr{$x_\mathrm{eq}^\sneg$} across \pr{$x_\mathrm{eq}^\uneg$},
approximated by Kramers' escape rate:
\begin{equation}
\kappa\approx \dfrac{\sqrt{\alpha\beta}}{2\pi}\exp\bigg(-\dfrac{\Delta U_0}{D}\bigg)\mbox{,}\label{eq:kramers}
\end{equation}
where $\alpha =
\partial_{yy}U_0(x_\mathrm{eq}^\sneg)$, $\beta =
-\partial_{yy}U_0(x_\mathrm{eq}^\uneg)$ and $\Delta U_r = U_r(y^\u[r])
- U_r(y^\s[r])$ ($\Delta U_0$ is then $\Delta U_r$ at $r=0$). The
approximation \eqref{eq:kramers} is accurate for $D\ll\Delta U_0$.

\subsection{First-order approximations of non-stationary Fokker-Planck
  Equation (FPE)}
\label{sec:approx}
\setcounter{paragraph}{0} We study the first-order deviation from the
quasi-stationary escape rate for $0<\epsilon\ll1$. For the following
section we fix a specific trajectory for $\mu(t)$ in $\dot
\mu=\epsilon\Gamma(\mu)$ by choosing the initial condition such that
\begin{displaymath}
  \mu=\mu_\mathrm{crit}\mbox{\quad for $t=0$}
\end{displaymath}
(remember that $\Gamma(\mu_\mathrm{crit})=1$ is the unique maximum of
$\Gamma$). We are interested in the parameter range where the diffusion
coefficient $D$ and the maximal ramp speed $r$ satisfy
\begin{align*}
  D&=O(\epsilon^{3/2})\mbox{,} & r-r_c(\epsilon)&=O(\epsilon)<0\mbox{.}
\end{align*}
Thus, the maximal ramp speed $r$ is close to, but below its critical
value $r_c(\epsilon)$ given in \eqref{eq:criticalrate}. \review{The
  scaling of the diffusion coefficient arises naturally from the
  change of coordinates to \eqref{ODE rc expansion} with white noise
  of variance $2D$; see Appendix \ref{app:crit rate expansion}. For $D
  \ll \epsilon^{3/2}$ we have a small-noise limit for $\epsilon \to 0$
  for the probability of escape, governed by Kramers' escape rate
  \eqref{eq:kramers}. For $D\gg\epsilon^{3/2}$ the probability of
  escape is dominated by noise-induced escape far away from the
  tipping. For $D\sim\epsilon^{3/2}$ all coefficients in the
  non-dimensionalized system are of order unity such that the
  time-dependence and the noise effects are in non-trivial balance for
  all small $\epsilon$}. The scaling implies in particular that
\begin{equation}
  \label{eq:rdrange}
  \pr{0<D\sim\Delta U_r\ll\Delta U_0\mbox{.}}
\end{equation}
\textbf{Remarks:} (\emph{a}) Condition~\eqref{eq:rdrange} \pr{means} that the escape rate is small
before and after the ramp ($t\to\pm\infty$, $\mu(t)$ far from
$\mu_\mathrm{crit}$) such that the Kramers approximation
\eqref{eq:kramers} for the escape rate is applicable for all times
outside of an interval $[t_0,T_{\mathrm{end}}]$ around $t=0$.

(\emph{b}) Condition~\eqref{eq:rdrange} also \pr{means} that the maximal ramp speed $r$
is sufficiently large such that approximation \eqref{eq:kramers} is no
longer true at the maximum speed of the ramp ($t=0$,
$\mu=\mu_\mathrm{crit}$), but it is still less than the saddle-node
rate $r_0$ (note that \pr{$\Delta U_{r_0}=0$}), the limit of the critical
rate $r_c(\epsilon)$ for $\epsilon\to0$. Thus, without noise ($D=0$),
there is a connection from $(x_\mathrm{eq}^\sneg,0)$ to
$(x_\mathrm{eq}^\spos,1)$ (the case of \emph{tracking} in
Section~\ref{sec:deterministic}). Let us pick one time profile $\tilde
x(t)$ on this connecting orbit.

\begin{figure}[h!]
        \centering       
                \includegraphics[scale = 0.3]{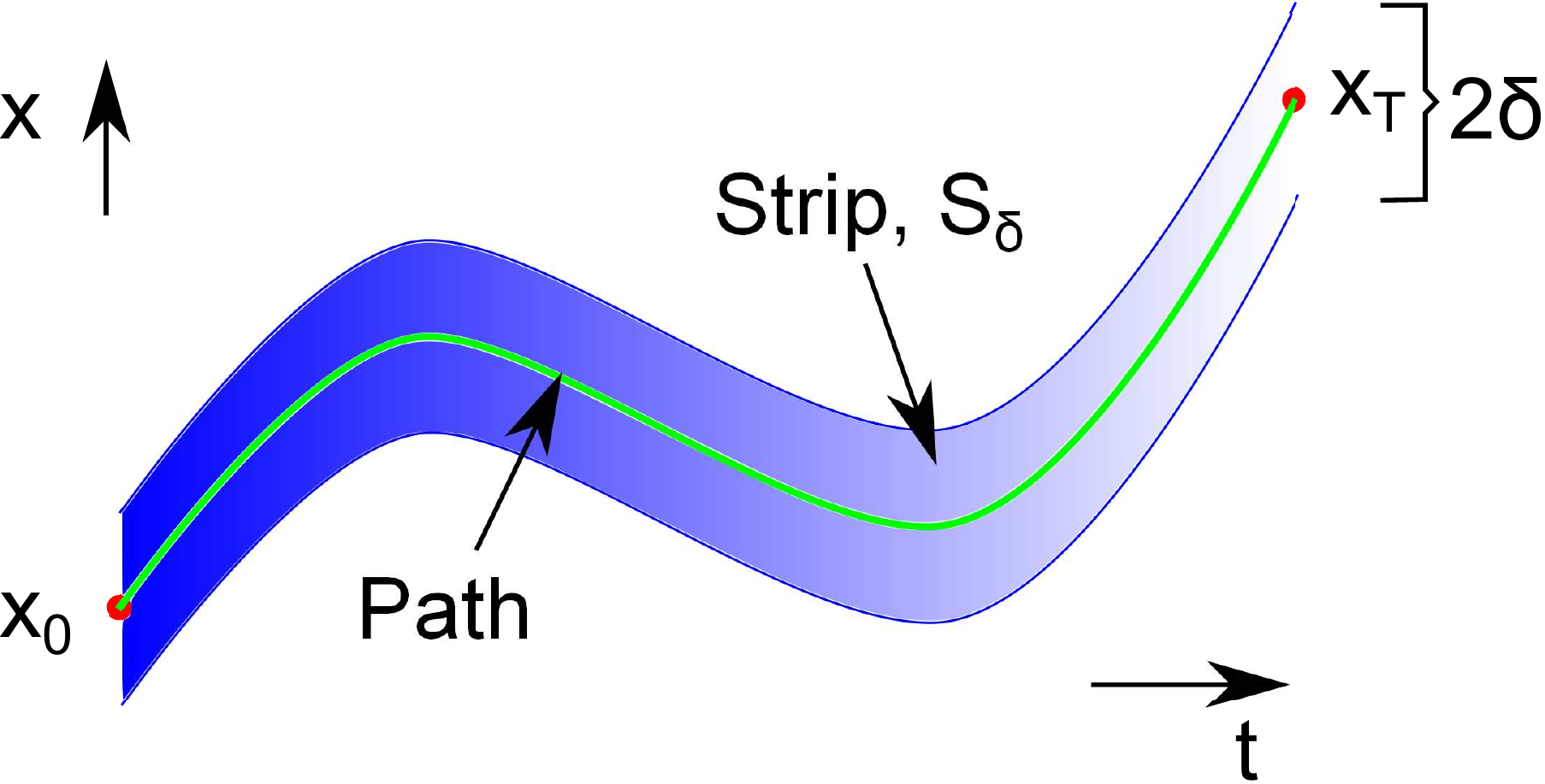}
                \caption{Sketch of the path $\tilde{x}(t)$ (green curve) and
                  the surrounding strip $S_\delta$ of width $2\delta$
                  (shaded blue). The probability density of a
                  realization passing through $(x,t)$ (always staying
                  within $S_\delta$) is $P(x,t)$.}
                \label{Sketch 2}
\end{figure}

We consider a starting position $(x,t) = (x_{0},t_{0})$ and an end
position at $(x,t) = (x_{T},T_{\mathrm{end}})$, and a strip
$S_\delta$ % The aim is to construct a method for calculating the probability of travelling along this path within a given band
of width $2\delta$ and length $T_{\mathrm{end}}-t_0$ around $\tilde x$,
see Figure~\ref{Sketch 2}.

% The method for calculating the probability of staying within the bands of the path is to use the probability density for a random variable $X_{t}$.
The Fokker-Planck equation
% , a linear parabolic partial differential equation (PDE), given as
\begin{equation}
\dfrac{\partial P}{\partial t} = D\dfrac{\partial^{2}P}{\partial x^{2}} - \dfrac{\partial}{\partial x}\bigg(f(x+r\mu(t)/\epsilon)P\bigg)\mbox{,}
\label{gen Fokker-Planck}
\end{equation}
describes the time evolution of the probability density $P(x,t)$ of
the random variable $X_{t}$, governed by \eqref{gen SDE}. If we impose Dirichlet boundary conditions,
\begin{align}
  0&=P(\tilde x(t)+\delta,t)\label{eq:dbcp}\mbox{,}\\
  0&=P(\tilde x(t)-\delta,t)\mbox{,} \label{eq:dbcm}
\end{align}
then $\int_a^e P(x,t)\d x$ is the probability that the solution of
\eqref{gen SDE}, starting at $t_0$ with probability
density $P(\cdot,t_0)$, is in $[a,e]$ at time $t$ and has never left
the strip $S_\delta$. % and evolve the density according to the Fokker-Planck
% equation \eqref{gen Fokker-Planck} on a domain corresponding to the
% band width with Dirichlet boundary conditions. The area of the density
% at $T_{\mathrm{end}}$ will reveal the probability for a realisation
% starting at $x_{0}$ to follow the path $\tilde{x}(t)$ within a given
% band of height $\delta$.

Consequently, the overall escape probability from the strip $S_\delta$ of width $2\delta$ around the path $\tilde x$ during time
interval $[t_0,T_{\mathrm{end}}]$ equals $1-\int_{-\delta}^\delta P(\tilde
x(T_{\mathrm{end}})+x,T_{\mathrm{end}})\d x$.  % However, we want to avoid solving a PDE and we can do
% this provided we assume that the probability density is
% quasi-stationary. If we first consider the stationary solution of the
% Fokker-Planck equation given by \citep{risken2012fokker}:

% \noindent where $P_{0}$ is a normalisation constant and $U(x,t) = -\int f(x,t)\,\mathrm{d}x$ is a potential.

\begin{figure}[h!]
        \centering
                {\includegraphics[scale = 0.4]{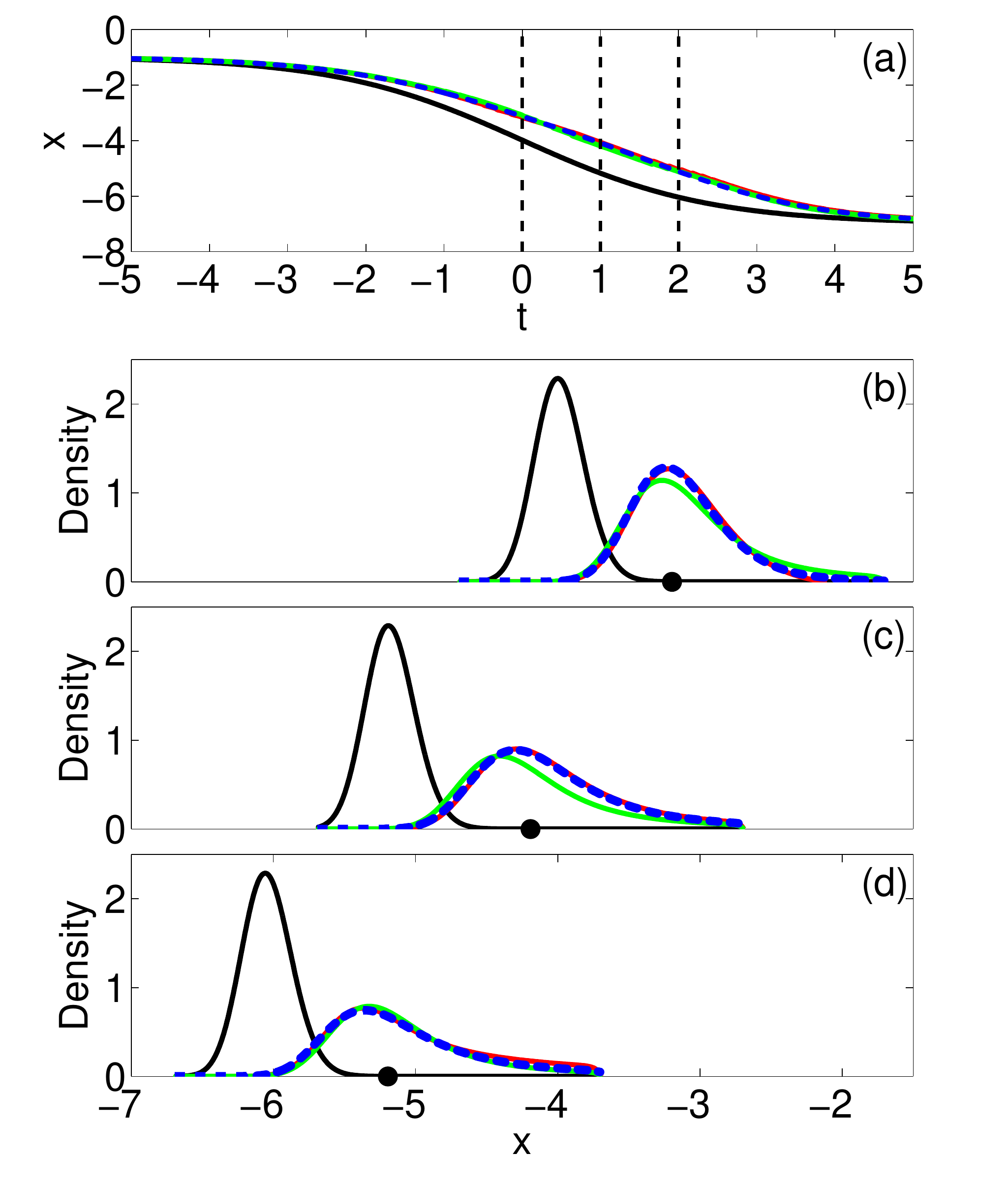}}
              \caption{Comparison of the
          single- (green (light gray)) and three-mode (red (dark gray)) approximations with the
          density from simulations (blue dashed) and stationary density
          (black) for \review{$\epsilon = 0.21$, $r = 1.26$}. Panel (a) provides the time profile for the location of the mean of each distribution. Vertical dashed lines indicate the times for which the densities are given in the remaining panels, namely (b) $t = 0$, (c) $t
          = 1$, and (d) $t = 2$. Black dot on $x$-axis (panels (b)-(d))
          corresponds to location of deterministic trajectory starting
          at $x_{0} = -1$ at $t = -10$. Parameters: Noise level $D =
          0.06$, width of strip $S_{\delta}$ $2\delta = 3$.} %$\rho = 0.2$ (a),(c),(e),(g) and
          %$\rho = 0.7$ (b),(d),(f),(h). Panels (a) and (b) give the time profile for the location of the mean of each distribution. Vertical dashed lines indicate the times for which the densities are given in the remaining panels, namely $t = 0$ - (c),(d), $t
          %= 1$ - (e),(f) and $t = 2$ - (g),(h). Black dot on $x$-axis (panels (c)-(h))
          %corresponds to location of deterministic trajectory starting
          %at $x_{0} = -1$ at $t = -10$. Parameters: Noise level $D =
          %0.06$, width of strip $S_{\delta}$ $2\delta = 3$.}
          \label{Density comparisons}
\end{figure}
Figure~\ref{Density comparisons} illustrates the shape of this
probability density $P(x,t)$ along the strip $S_\delta$ for a ramp
speed $r < r_{c}(\epsilon)$. The figure uses parameters from the
specific example introduced in Section~\ref{sec: Setup}.
 The first moment of $P(x,t)$ (the mean)
is shown in panel (a) of Figure~\ref{Density comparisons}. Panels
(b)--(d) show the profile of $P(x,t)$ for selected times $t$.
%two different different ramp speeds $r$ and
%perturbation parameters $\epsilon$ (keeping $r/\epsilon$ constant):
%$\epsilon\ll1$ in the left column, and $\epsilon$ larger (but $r$
%still less than $r_c(\epsilon)$) in the right column. 
The numerical solution of \eqref{gen Fokker-Planck}--\eqref{eq:dbcp}
is shown as a blue (dashed) curve. The other curves are the
approximations described below.
\paragraph{Uncentered quasi-stationary density}
\label{sec:uncentered}
The crudest approximation assumes that the density is approximately
stationary throughout the ramp. This implies that $\partial_tP$ is
small in \eqref{gen Fokker-Planck}. Replacing $\partial_tP$ with zero
in \eqref{gen Fokker-Planck} and imposing a Dirichlet boundary
condition on the right end ($P(\tilde x(t)+\delta,t)=0$), the solution
$P$ of \eqref{gen Fokker-Planck} has the form
\begin{equation}\label{eq:pstar}
  P_*(x,t) = P_{0}(t)\int\limits_x^{\tilde x(t)+\delta}
  \exp\left[\frac{U(x',t)-U(x,t)}{D}\right]\d x'\mbox{,}
\end{equation}
where $\partial_xU(x,t)=-f(x+r\mu(t)/\epsilon)$ for each fixed $t$.
% $P_{*}$ is a solution of the stationary form of \eqref{gen
%   Fokker-Planck} (setting $\partial_tP=0$), \js{imposing a Dirichlet
%   boundary condition only on the right side: $P(\tilde
%   x(t)+\delta,t)=0$}.
The spatial shape of $P_*(\cdot,t)$ is nearly unchanged, only shifted
by $r\mu(t)/\epsilon$ for different times $t$. As the density has
nearly constant shape for all $t$, the escape rate is nearly constant
in time as well (hence, it is equal to $\kappa$ by remark
(\emph{a})). This escape rate determines the normalization
  constant $P_0(t)$: $P_0(t)\approx P_\infty(1-(t-t_0)\kappa)$, where
  $P_\infty$ is such that the initial density \pr{$P_{*}(x,t_0)$} has
  a unit integral. This approximation, shown in black in
  Figure~\ref{Density comparisons} for the example from
  Section~\ref{sec: Setup}, does not catch the effect of a non-zero
  $r$: it is (nearly) independent of $r$ (becoming independent of $r$
  in the limit $\delta\to\infty$). The density \pr{$P_{*}(x,t)$} is
  centered at \pr{$x_\mathrm{eq}^\sneg-r\mu(t)/\epsilon$}, which is
  visibly smaller than $\tilde x(t)$ (the location of the
  deterministic trajectory, highlighted by the black dot on the
  $x$-axis in \pr{panels (b)--(d)} of Figure~\ref{Density
    comparisons}).

\paragraph{Instantaneous eigenmodes of FPE, centered at $\tilde x$}
%  approximations to the probability
% density $P(x,t)$ from the numerical solution of \eqref{gen
%   Fokker-Planck} (blue) at different times (rows) and different ramp
% speeds $\epsilon$ (columns) for the specific example introduced in
% Section~\ref{sec: Setup}.  For a slow drift speed (left column) the
% numerical solution in blue has a slightly wider distribution and also
% lags the stationary density (solving \eqref{gen Fokker-Planck} with
% $\partial_tP=0$ and boundary condition \eqref{eq:dbcp} at the upper
% boundary only), shown in black, This is also highlighted by the black
% dot on the $x$-axis of the figures, which denotes the location of the
% deterministic trajectory $\tilde x(t)$. These features are exaggerated
% further for a larger drift speed (right column).  We therefore
% conclude that the stationary density is poor in approximating the true
% quasi-stationary probability density. Instead w

The instantaneous eigenmode expansion for the linear operator of the
Fokker-Planck equation \eqref{gen Fokker-Planck} follows an approach
similar to that presented in
\citet{risken2012fokker,zhang1997distributed} but for a time dependent
deterministic part $f(x+r\mu(t)/\epsilon)$ instead of a
time-independent $f(x)$.  Figure~\ref{Density comparisons} shows the single-mode approximation ($n=1$) in green (light gray) and the three-mode
approximation ($n=3$) of the probability density in red (dark gray).
% \change{At this stage we will simply highlight that
%   the mode approximations provide a better approximation to the
%   probability density than the stationary density. Further analysis
%   will be considered later for the specific example they apply to.}

We first change to a co-moving coordinate system $x=\tilde x(t)+y$ with respect to $y$, such that, the SDE~\eqref{gen SDE} has the form
\begin{align}
  \label{eq:shiftSDE}
  \d y&=-y g(y,t) \d t +\sqrt{2D}\d W_t\mbox{, where\ }\\
  g(y,t)&=-\int_0^1f'(sy+\tilde y(t))\d s\nonumber\mbox{ with\ }\\
  \tilde y(t)&=\pr{y^\s[r\Gamma(\mu(t))]}-\epsilon
  \frac{r\Gamma'(\mu(t))\Gamma(\mu(t))}{\pr{f'(y^\s[r\Gamma(\mu(t))])^2}}+O(\epsilon^2)\mbox{}\nonumber
\end{align}
(recall that $y^\s[r\Gamma(\mu)]$ is the stable equilibrium of $\dot
y=f(y)+r\Gamma(\mu)$ with fixed $\mu$). For sufficiently small
$\epsilon$ and $r<r_0$ we have that in the moving coordinates $y$ the
equilibrium $y=0$ is stable without noise for each fixed $t$ (by the
stability assumption on \pr{$y^\s[r\Gamma(\mu(t))]$} and because
$r<r_0$). That is, for $t\in[t_0,T_{\mathrm{end}}]$
\begin{displaymath}
  g(0,t)=-f'(\tilde y(t))=-f'(y^\s[r\Gamma(\mu(t))])+O(\epsilon)>0\mbox{.}
\end{displaymath}
The Fokker-Planck equation for the density $P(y,t)$ over $y$ operates
then on the fixed domain $[-\delta,\delta]$:
\begin{equation}\label{eq:shiftFP}
 \pr{\dfrac{\partial P}{\partial t}= D\dfrac{\partial^{2}P}{\partial y^{2}} + \dfrac{\partial}{\partial y}\bigg(yg(y,t))P\bigg)=: A(t)P}
\end{equation}
with Dirichlet boundary conditions $P(-\delta,t)=P(\delta,t)=0$. The
operator $A(t)$ is self-adjoint with respect to the scalar product
\begin{equation}\label{eq:scalarproduct}
  \langle w,v\rangle_t=\int_{-\delta}^\delta w(y)v(y)
  \exp\left(\frac{U(y,t)}{D}\right)\,\d y\mbox{,}
\end{equation}
where $U$ is the potential corresponding to the drift $-yg(y,t)$
($\partial_yU(y,t)=yg(y,t)$). Since $g(y,t)>0$, this effective
potential $U(y,t)$ has a \pr{critical point} at $y=0$. For each fixed $t\in[t_0,T_{\mathrm{end}}]$,
the spectrum of $A(t)$, shown \pr{later in Figure~\ref{Eigenvalues all t} for the specific
example considered} in Section~\ref{sec: Setup}, consists of eigenvalues
$\gamma_i(t)$ with eigenfunctions $v_i(y,t)$:
\begin{equation}\label{eq:FPeig}
   \gamma_i v_i = D\dfrac{\partial^{2}v_i}{\partial y^{2}} + \dfrac{\partial}{\partial y}\left[yg(y,t)v_i\right]=A(t)v_i\mbox{.}
\end{equation}

% We ascertain the instantaneous eigenmodes by first
% simplifying the Fokker-Planck equation \eqref{gen Fokker-Planck} down
% to an ordinary differential equation (ODE). This is using the
% assumption that the probability density is quasi-stationary, therefore
% we represent the Fokker-Planck equation as an eigenvalue parameter
% dependent problem:

% \begin{equation*}
% \gamma(t) P(x,t) = \bigg[D\dfrac{\partial^{2}}{\partial x^{2}} - \dfrac{\partial f(x,t)}{\partial x} - f(x,t)\dfrac{\partial}{\partial x}\bigg]P(x,t) 
% \end{equation*}

% \noindent we see that this can be expressed in terms of an unbounded linear operator $A(t)$ and so we have

% \begin{equation}
% \dfrac{\partial P(x,t)}{\partial t} = \gamma(t) P(x,t) = A(t)P(x,t)
% \label{eigenvalue Fokker-Planck}
% \end{equation}

% \noindent where

% \begin{equation}
% A(t) = D\dfrac{\partial^{2}}{\partial x^{2}} - \dfrac{\partial f(x,t)}{\partial x} - f(x,t)\dfrac{\partial}{\partial x}
% \label{linear operator}
% \end{equation}

% \noindent The eigenvalues $\gamma_{i}(t)$ of the linear operator $A(t)$ and the corresponding instantaneous eigenmodes $v_{i}(t)$ satisfy the following relation \citep{peube2009fundamentals}:

% \begin{equation*}
% A(t)v_{i}(t) = \gamma_{i}(t)v_{i}(t)
% \end{equation*}

\noindent The eigenfunctions $v_{i}(\cdot,t)$ (called
\emph{instantaneous} modes as they are time-dependent) form an
orthonormal basis of $L^{2}$ with respect to
$\langle\cdot,\cdot\rangle_t$. Thus, we can expand the solution
$P(y,t)$ of \pr{\eqref{eq:shiftFP}} as a linear combination
of the instantaneous eigenmodes $v_i$:
\begin{equation}
P(y,t) = \sum_{i = 1}^{\infty}a_{i}(t)v_{i}(y,t)\mbox{,}
\label{Assumed solution}
\end{equation} 
where the $a_{i}(t)$ are scalar coefficients at each time $t$.
%  and we impose Dirichlet boundary conditions:

% \begin{equation*}
% P(\tilde{x}(t)-x_{l},t) = P(\tilde{x}(t)+x_{u},t) = 0 \qquad \forall t
% \end{equation*}

% \noindent where $x_{l}$, $x_{u}>0$ determine the location of the band relative to the path $\tilde{x}(t)$. 
% The initial values $a_{i}(0)$ from equation \eqref{Assumed solution}
% are given by projection of the given initial density $P(y,0)$
% \citep{williams2012linear} using the scalar product $\langle\cdot,\cdot\rangle_t$ at $t=0$: 
% \begin{displaymath}
% a_{i}(0) = \langle v_{i}(\cdot,0),P(\cdot,0)\rangle_0\mbox{.}
% \end{displaymath}
% \noindent For this the eigenmodes $\textbf{w}_{k}$ of the adjoint of the linear operator $A$ are required. The adjoint operator $A^{adj}$ is given by:

% \begin{equation*}
% A^{adj}(t) = D\dfrac{\partial^{2}}{\partial x^{2}} + f(x,t)\dfrac{\partial}{\partial x}
% \end{equation*} 
Inserting the expansion \eqref{Assumed solution} into equation \pr{\eqref{eq:shiftFP}}, applying $\langle v_k,\cdot\rangle_t$, and truncating at a finite $n$ gives \citep{williams2012linear}:
% \begin{equation}
% \sum_{i = 1}^{\infty}(\dot{a}_{i}\textbf{v}_{i} + a_{i}\dot{\textbf{v}}_{i}) = \sum_{i = 1}^{\infty}\gamma_{i}(t)a_{i}(t)\textbf{v}_{i}(t)
% \label{calculation of akdot}
% \end{equation}

% \noindent Multiply equation (\ref{calculation of akdot}) on the left by $\textbf{w}_{k}^{T}$ and use the definition the eigenmodes of $A$ and $A^{adj}$ are orthonormal to each other i.e. $\textbf{w}_{i}^{T}\textbf{v}_{j} = \delta_{i,j}$:
\begin{align}
\dot{a}_{k} &= \gamma_{k}(t)a_{k}(t) - \sum_{i = 1}^n\langle v_k(t),\dot v_i(t)\rangle_t a_{i}(t)\mbox{,}
\label{akdot}\\
a_{i}(0) &= \langle v_{i}(\cdot,0),P(\cdot,0)\rangle_0\mbox{,}\nonumber
\end{align} 
where the coupling coefficients $\langle v_k(t),\dot v_i(t)\rangle_t$
are of order $\epsilon/|\gamma_k(t)-\gamma_i(t)|$ for $i\neq k$. Thus, the sum is
convergent and the truncated solution 
\begin{equation}\label{eq:Pn}
P_{n}(x,t) = \sum_{i = 1}^{n} a_{i}(t)v_{i}(x,t)
\end{equation}
converges to $P$ for $n\to\infty$ and for $\epsilon\to0$
($-\gamma_k=O(k^2)$ for positive $D$ \pr{(see Figure~\ref{Eigenvalue spectrum all t} for the specific example)}, so, in particular, $P_n-P=O(\epsilon)$).
% is because the
% $\dot{\textbf{v}}_{i}(t)$ are close to zero if there is very little
% change in the eigenmodes from one time step to the next.

% We truncate the solution in (\ref{Assumed solution}) into two parts because \change{the scalars $a_{k}(t)$ in front of the latter eigenmodes are small:}

% \begin{equation*}
% P(x,t) = P_{n}(x,t) + P_{r}(x,t)
% \end{equation*}

% \noindent where 

% \begin{equation*}
% P_{n}(x,t) = \sum_{i = 1}^{n} a_{i}(t)v_{i}(x,t)
% \end{equation*}

% \noindent contains the sum of the first $n$ eigenmodes and $P_{r}(x,t)$ the sum of the remaining eigenmodes. In the limit we have:

% \begin{equation*}
% |P_{n}(x,t) - P(x,t)| \rightarrow 0 \qquad \text{as} \quad n \rightarrow \infty \qquad \forall t>0
% \end{equation*} 

Since the initial $t_0$ is such that $\mu(t_0)$ is still close to
$0$, $\gamma_1(t_0)$ will be very close to zero by \pr{remark}
(\emph{a}) that escape is unlikely outside of the ramping time
interval. \pr{Remark} (\emph{a}) also implies that the coefficients of
the initial value $P(\cdot,t_0)$, $a_i(t_0)$, are close to zero for
$i>1$. % Furthermore, for small $\epsilon>0$ the effect of the coupling
% terms on the solution is small such that truncation at small $n$ is
% justified for small $\epsilon$.

Figure~\ref{Density comparisons} illustrates that the truncation error
for small $n$ occurs in the tails of the distribution. For example, for the
single-mode approximation with $n=1$ the equation \eqref{akdot}
simplifies to
\begin{displaymath}
  \dot a_1=[\gamma_1-\langle v_1,\dot v_1\rangle_t]a_1\mbox{.}
\end{displaymath}
Ignoring the $O(\epsilon)$ term $\langle v_1,\dot v_1\rangle_t$, the
single-mode approximation results in an approximate solution
\begin{displaymath}
\pr{P_1(y,t)=\exp\bigg(\int_{t_0}^{T_{\mathrm{end}}}\gamma_1\mathrm{d}t\bigg)v_1(y,t)}.
\end{displaymath}
Thus, truncation at $n=1$
assumes that the density instantaneously adjusts its shape to the
shape of the effective potential well \pr{$U(y,t)=-\int yg(y,t)\mathrm{d}y$} at every
time $t$.

For each particular truncation $n$, the probability
$\mathbb{P}_{M}$ for the trajectory of a realization to \prnew{not} remain within
the strip $S_\delta$ is approximately
\begin{equation}
\mathbb{P}_{M} = \prnew{1 - }\int_{-\delta}^\delta P_n(y,T_{\mathrm{end}})\,\mathrm{d}y\mbox{.}
\label{Modes prob}
\end{equation}

\noindent
\subsection{Perturbation approximation of the dominant eigenvalue}
The dominant eigenvalue $\gamma_{1}(t)$ and eigenfunction $v_1$ can be
approximated via a linear perturbation analysis from the small-noise
limit ($D\to0$). Hence, we can approximate the dominant term in
$\mathbb{P}_{M}$ for the truncation $n=1$, which is accurate to order
$\epsilon$. Consider again the eigenvalue problem for the
Fokker-Planck equation:
\begin{equation}
\gamma(t) P(y,t) = D\dfrac{\partial^{2}P(y,t)}{\partial y^{2}} + \dfrac{\partial }{\partial y}[U'(y,t)P(y,t)]\mbox{,}
\label{eigenvalue FP}
\end{equation}
for $y\in(a,\delta)$ (where $a$ can be $-\infty$ in some
expressions below, but we will finally set $a=-\delta$). We now
consider $t$ simply as a parameter in the eigenvalue problem (such
that eigenvalue $\gamma$ and eigenfunction $P$ depend on the parameter
$t$ since the coefficient $U'$ depends on $t$). We will drop this
parameter $t$ throughout this subsection. The basic building block of
solutions of \eqref{eigenvalue FP} is the function $\exp(-U(y)/D)$,
which we call
\begin{equation}\label{eq:pdef}
  p(y)=\exp(-U(y)/D)\mbox{,}
\end{equation}
along with anti-derivatives of products of $p$ of various orders,
which we call
\begin{displaymath}
  p_{s_1\ldots s_k}(y_0)=\int\limits_{y_0}^\delta\!\!\ldots\!\!\!\!\int\limits_{y_{k-1}}^\delta\!\!
  \exp\left[\sum_{j=1}^k\frac{(-1)^{s_j}U(y_j)}{D}\right]\d y_k\ldots\d y_1
\end{displaymath}
(the subscripts $s_j$ will be $1$ or $2$). We know that
\begin{equation}
  \gamma=0\mbox{,}\quad P_*(y) = p(y)p_2(y)/p_{12}(a)\mbox{}
\label{initial sol of EFP}
\end{equation}  
solve \eqref{eigenvalue FP} with the two
boundary/integral conditions
\begin{equation}
P(\delta) = 0, \qquad \int_{a}^{\delta} P(y)\mathrm{d}y = 1\mbox{.}
\label{eigenvalue problem conds stat}
\end{equation} 
The expression for $P_*$ in \eqref{initial sol of EFP}, equals expression
\eqref{eq:pstar} for $P_*$ with the specific normalization constant
\prnew{$P_0=J/D$ where $J$ represents the probability flux. The probability flux $J$ is the flow of probability per unit time per unit area.
%$P_0=1/p_{12}(a)$.
%% Equation~\eqref{eigenvalue FP} with \eqref{eigenvalue problem conds
%%   stat} has the solution
%In the small-noise limit, the anti-derivatives of $p$ can be
%approximated further in \eqref{initial sol of EFP}, leading to the
%well-known approximation for Kramers' escape rate.

  The integral condition in \eqref{eigenvalue problem conds stat} is
  based on an assumption that is only approximately correct if the
  noise level $D$ is small: the probability flux $J$ is constant in $y$
  such that the flux through the right boundary at $+\delta$ must also
  enter at $y=-\infty$ such that $J$ is given by
\begin{equation}
J=\frac{D}{p_{12}(a)}.
\label{probability flux}
\end{equation}} % If one wants to approximate a stationary distribution with
  % an ensemble in a Monte-Carlo simulation of the SDE \eqref{gen SDE},
  % one has to re-initialize those realizations that escaped. The
  % conditions \eqref{eigenvalue problem conds stat} correspond to a
  % re-initialization of the realization at $x=-\infty$,
  The more appropriate boundary conditions for non-small noise level
  $D$ result in % there-initialization is at a random point in $(-\delta,\delta)$
  % according to the current density $P$. The difference between those
  % re-initializations is of higher order only in the limit of
  % $D\to0$. For non-small $D$, the solution of
  the eigenvalue problem \eqref{eigenvalue FP} for $\gamma$ and $P$
with
\begin{equation*}
P(a)=0,\qquad P(\delta) = 0, \qquad 
\int_{a}^{\delta} P(y)\mathrm{d}y = 1\mbox{,}
\end{equation*} 
which leads to a uniformly non-zero $\gamma$, including for the limit
$a\to-\infty$.

We can express a first-order approximation of $\gamma$ for non-small
noise in terms of $p$, given in \eqref{eq:pdef}, by treating it as a
perturbation of the small-noise limit and of $\gamma=0$, $P=P_*$. For
a finite $a\ll-1$ let us introduce the value of the solution $P$ of
\eqref{eigenvalue FP} at $a$ as a parameter $\pi_a$:
\begin{equation}\label{eq:bca}
  P(a)=\pi_a\mbox{.}
\end{equation}
Then, we get a solution pair $(\gamma,P)$ of \eqref{eigenvalue FP}
with boundary conditions \eqref{eigenvalue problem conds stat},
\eqref{eq:bca} for each small $\pi_a$. For
$\pi_a=P_*(a)=p(a)p_2(a)/p_{12}(a)$, the solution is $\gamma=0$,
$P=P_*$. So, to first order in $\pi_a$, we have for the parameter
$\pi_a=0$
\begin{equation}\label{eq:gamma:1st}
  \gamma\approx-\gamma'P_*(a)=-\gamma'\frac{p(a)p_2(a)}{p_{12}(a)}\mbox{\ with\ }
\gamma'=\frac{\d\gamma}{\d \pi_a}\biggl\vert_{\pi_a=P_*(a)}\mbox{.}
\end{equation}
The scalar $\gamma'$ is part of the solution pair $(\gamma',q)$ of the
linearization of eigenvalue problem \eqref{eigenvalue FP} with
boundary conditions \eqref{eigenvalue problem conds stat},
\eqref{eq:bca}, with respect to $\pi_a$ in $\pi_a=P_*(a)$, $P=P_*$,
$\gamma=0$:
\begin{equation}
\gamma'P_{*} = \dfrac{\partial}{\partial y}\bigg(D\dfrac{\partial q}{\partial y} + U'(y)q\bigg)
\label{lin eig eq}
\end{equation}
with conditions
\begin{equation*}
q(a) = 1, \qquad q(\delta) = 0, \qquad \int_{a}^{\delta} q(y)\mathrm{d}y = 0.
\end{equation*}
This is an affine equation for $\gamma'$ and $q$, which can be solved
by integration, resulting in
\begin{align*}
  \gamma'=\frac{Dp_{12}(a)^2}{p(a)[p_{12}(a)p_{212}(a)-p_2(a)p_{1212}(a)]}\mbox{,}
\end{align*}
such that the first-order estimate \eqref{eq:gamma:1st} gives
\begin{equation}\label{lead eig}
  \gamma_1:=\gamma\approx
  \frac{D}{\displaystyle\frac{p_{1212}(a)}{p_{12}(a)}-\frac{p_{212}(a)}{p_2(a)}}\mbox{.}
\end{equation}
Taking into account now that all quantities in \eqref{lead eig} depend parametrically on time $t$, the probability $\mathbb{P}_{P}$ of not following a path
within a specific region is approximately 
\begin{equation}
\mathbb{P}_{P} = \prnew{1 -} \exp\bigg(\int_{t_0}^{T_{\mathrm{end}}} \gamma_{1}(t)\d t\bigg)\mbox{,}
%\exp\bigg(\mathrm{d}t\sum\limits_{k=0}^{N+1} \gamma_{1}(t_{k})\bigg)
\label{eig prob}
\end{equation}
where $\gamma_1$ is given approximately in \eqref{lead eig}, when
inserting $-\delta$ for $a$. We can compare \eqref{eig prob} with the simpler probability formula $\mathbb{P}_{J}$ which is valid in the small noise limit
\begin{equation}
\mathbb{P}_{J} = 1 - \exp\bigg(-\int_{t_{0}}^{T_{\mathrm{end}}}J(t)\,\mathrm{d}t\bigg)\mbox{,}
\label{Flux prob}
\end{equation}
where the probability flux $J$ is given in \eqref{probability flux}, when again inserting $-\delta$ for $a$.
% \noindent where the time interval has been split into $N+1$ grid points $t_{k}$ all separated by a distance $\mathrm{d}t$.

\pr{Section \ref{sec: Probability}} will compare Monte-Carlo simulations, the
numerical approximation using the first $n$ instantaneous eigenmodes
of the linear operator $A(t)$ of the Fokker-Planck equation, the
perturbation formula for the leading eigenvalue $\gamma_{1}(t)$, and the
formula for the probability flux $J=D/p_{12}(a)$ for the stationary
density $P_*$ in \eqref{initial sol of EFP}, which is accurate for
small escape rates.

%  The example will be the prototype

% We will now evaluate these approximations for a specific example from noise and rate-induced tipping.

\section{Saddle-node normal form with parameter ramp and noise}
\label{sec: Setup}

A prototypical model for rate-induced tipping was introduced by
\citet{ashwin2012tipping}. % We modify the notation here to adapt it to
% the convention used in sections~\ref{sec:deterministic} and
% \ref{sec: Methods}.
The time evolution of a scalar dependent variable
$x(t) \in \mathbb{R}$ is described by the saddle-node normal form
equation:

\begin{equation}
\dot{x} = f(x,\lambda) = (x+\lambda)^{2} - 1\mbox{,}
\label{saddle node nf}
\end{equation}

\noindent where w.l.o.g. we have set the normal form parameter to
equal $1$. The ODE \eqref{saddle node nf} has two families of
equilibria; one stable family \pr{$W^\s[0] = -\lambda -1$} and one
unstable family \pr{$W^\u[0] = -\lambda + 1$}. The parameter $\lambda$
in equation \eqref{saddle node nf} is assumed to be time dependent
following a ramp given by:

\begin{equation}
\lambda(t) = \dfrac{\lambda_{\max}}{2}\bigg(\tanh\bigg(\dfrac{\lambda_{\max}\rho t}{2}\bigg) + 1\bigg)\mbox{,}
\label{ramp eq}
\end{equation}   

\noindent where $\lambda_{\max}$ determines how far the parameter $\lambda$ is
shifted and $\rho$ adjusts the speed of the ramp. Equation
\eqref{ramp eq} can be described by an ODE for $\lambda$ with the
condition $\lambda(0) = \lambda_{\max}/2$. Therefore the prototypical
model can be described by the two dimensional ODE in the $(x,\lambda)$
phase plane

\begin{align}
\label{sn ODE}
\dot{x} &= (x+\lambda)^{2} - 1\mbox{,} \\
\dot{\lambda} &= \rho\lambda(\lambda_{\max} - \lambda)\mbox{.}
\label{lambda ODE}
\end{align}

\noindent For this system of ODEs a critical speed $\rho =
\rho_{c}=4/[\lambda_{\max}(\lambda_{\max}-2)]$ was found in
\citep{perryman2015tipping}, at which a heteroclinic connection
$(x,\lambda)=(-1+(2/\lambda_{\max}-1)\lambda,\lambda)$ from $(-1,0)$ to \pr{$(1-\lambda_{\max},\lambda_{\max})$}
occurs (setting the critical rate $\rho_c$ for rate-induced
tipping). The time profile and phase portrait for
$\rho<\rho_{c}$ is presented in Figure~\ref{Deterministic
  plots}. For a complete overview of all possible time profiles and
phase portraits, see \citet{ritchie2015early}.

\begin{figure}[ht]
        \centering
        %\begin{subfigure}[h!]{0.45\textwidth}
        \subcaptionbox{\label{TP eps small}}[0.45\linewidth]
                {\includegraphics[scale = 0.3]{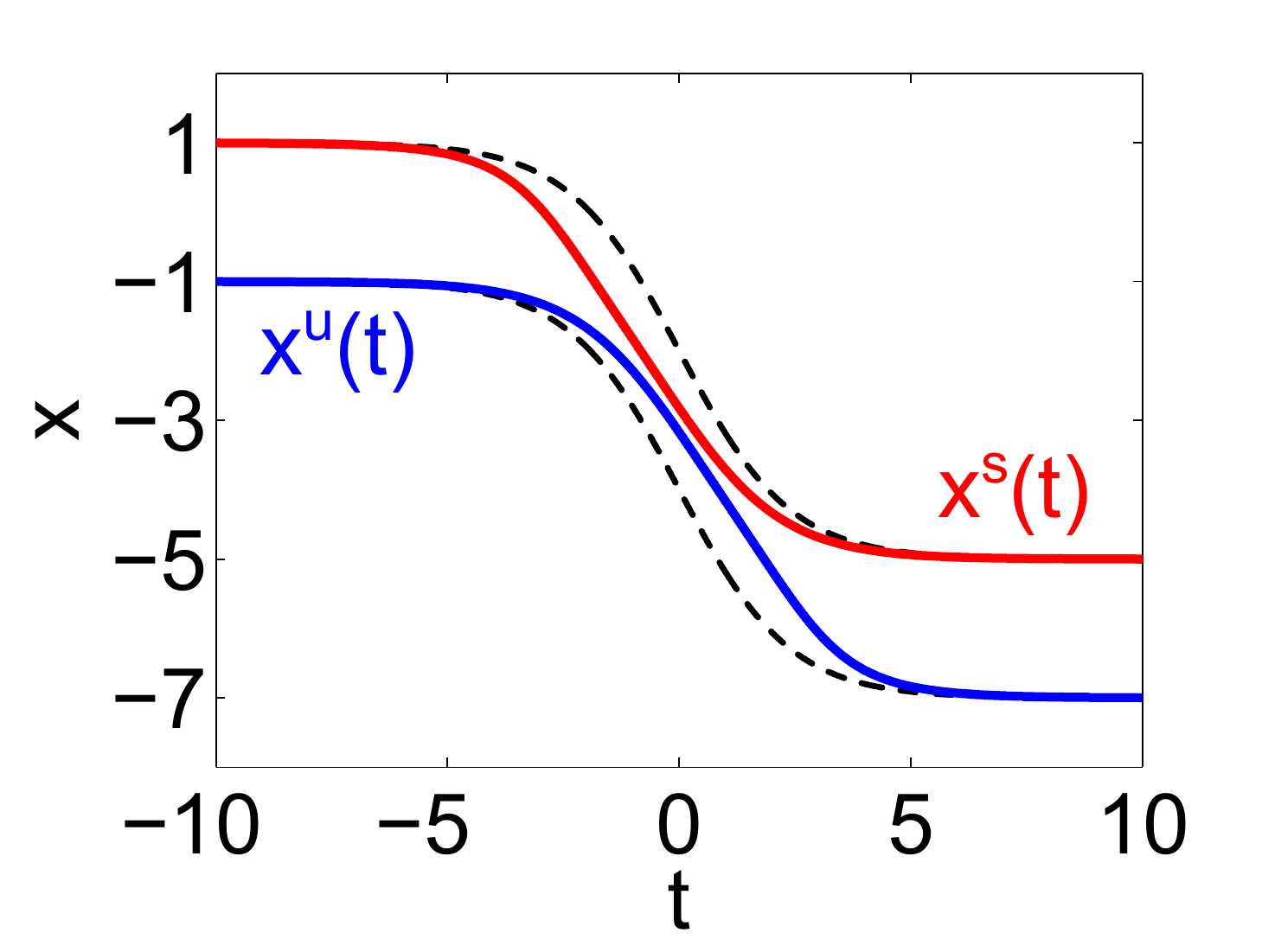}}
                %\caption{Time profile of manifolds $W^{u}(S_{-})$ and $W^{s}(U_{+})$}
                %\label{TP eps small}
        %\end{subfigure}%
       \hfill %add desired spacing between images, e. g. ~, \quad, \qquad, \hfill etc.
        %\begin{subfigure}[h!]{0.45\textwidth}
        \subcaptionbox{\label{PP eps small}}[0.45\linewidth]
                {\includegraphics[scale = 0.3]{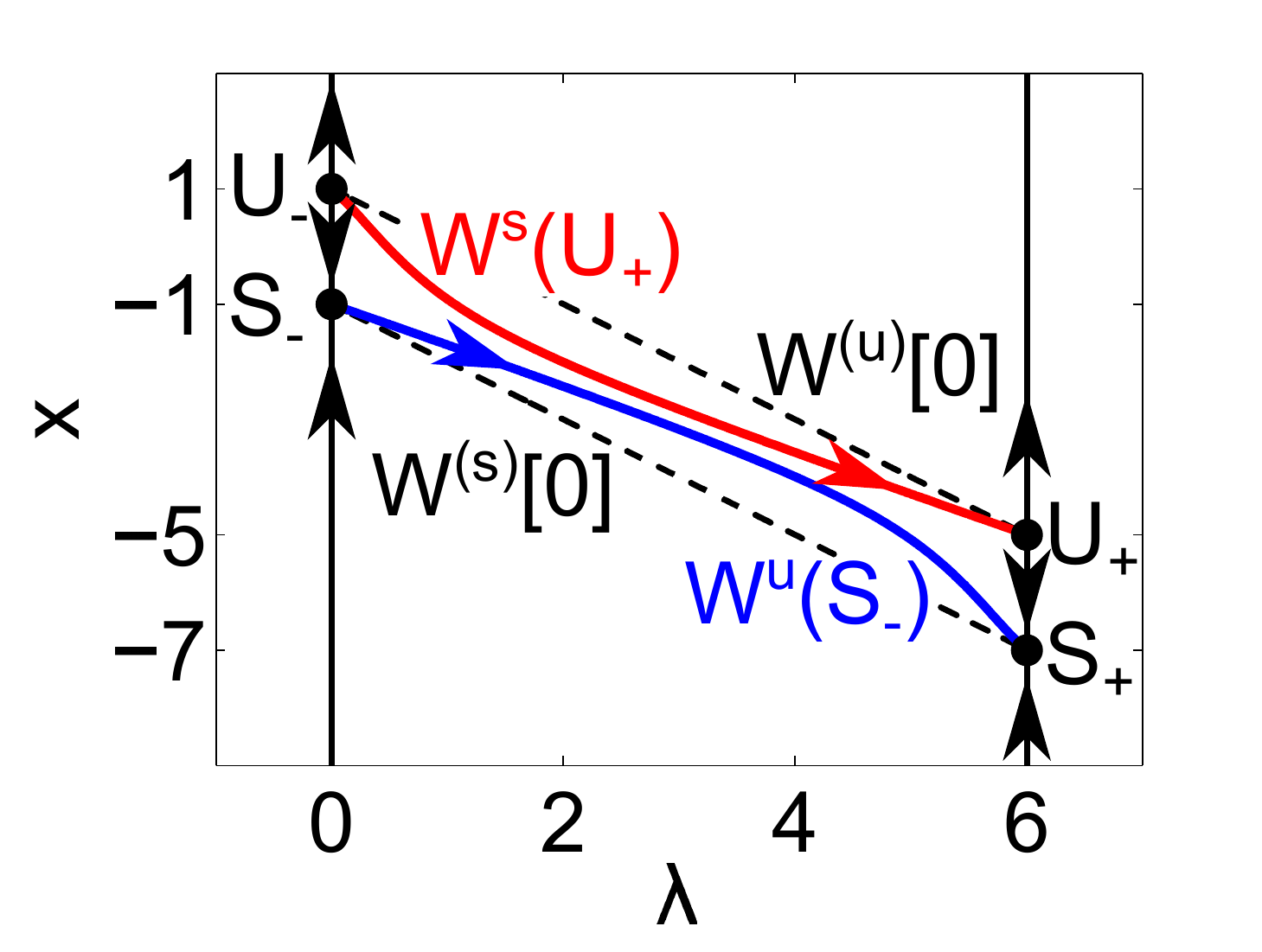}}
                %\caption{Phase plane of system (\ref{rtip ODE}),(\ref{rtip lambda dot init})}
                %\label{PP eps small}
        %\end{subfigure}
        ~ %add desired spacing between images, e. g. ~, \quad, \qquad, \hfill etc.
        \caption{Time profile (a) and phase plane
            (b) of system \eqref{sn ODE}--\eqref{lambda ODE} for \pr{$\rho
            = 0.14 < \rho_{c}$}. Black dashed curves are the stable
            \pr{($W^{(s)}[0]=-1-\lambda$) and unstable ($W^{(u)}[0]=1-\lambda$) branches of} equilibria in the
            limit $\rho = 0$, blue and red curves are the unstable
            and stable manifolds, $W^{u}(S_{-})$ and $W^{s}(U_{+})$,
            respectively \pr{($\lambda_{\max}=6$)}.}\label{Deterministic plots}
\end{figure} 

System \eqref{sn ODE}--\eqref{lambda ODE} has 4 equilibria: two
saddles $S_{-}=(-1,0)$, \pr{$U_{+}=(-5,6)$}, one stable node
\pr{$S_{+}=(-7,6)$} and one unstable node $U_{-}=(1,0)$; see
Figure~\ref{PP eps small}. The dashed lines $W^{(s)}[0]=-1-\lambda$
and $W^{(u)}[0]=1-\lambda$ represent the family of stable and
unstable equilibria for $\rho = 0$ respectively. The curve
$W^{u}(S_{-})$ is the unstable manifold of the saddle $S_{-}$ and
$W^{s}(U_{+})$ is the stable manifold of the saddle $U_{+}$. The time
profile for $x$ on the invariant manifolds $W^{u}(S_{-})$ and
$W^{s}(U_{+})$, denoted $x^{u}(t)$ and $x^{s}(t)$ respectively, is
given in Figure~\ref{TP eps small}.

The manifold $W^{s}(U_{+})$ acts as a separatrix partitioning the plane into two distinct regions. Below $W^{s}(U_{+})$ all trajectories are attracted towards the stable node $S_{+}$, while any trajectories above the separatrix escape to $+\infty$ in finite time. 

\pr{We will continue using the parameters $\rho$ and $\lambda_{\max}$ to remain consistent with the studies of \citep{ashwin2012tipping,perryman2015tipping,ritchie2015early}.} The parameters $\rho$ and $\lambda_{\max}$ (both are non-small) have the following
relation to the parameters $\epsilon$ and $r$ of
the sections \ref{sec:deterministic} and \ref{sec: Methods}:
\begin{align*}
  \epsilon&=\frac{\rho\lambda_{\max}}{4}\mbox{,} &
  r&=\frac{\rho\lambda_{\max}^2}{4}\mbox{,} &
  \lambda_{\max}&=\frac{r}{\epsilon}\mbox{,} &
  \rho&=\frac{4\epsilon^2}{r}\mbox{.}
\end{align*}
Using $r$ and $\epsilon$, system~\eqref{sn ODE}--\eqref{lambda ODE} has
the form
\begin{align*}
  \dot x&=(x+r\mu/\epsilon)^2-1\mbox{,} &
  \dot \mu&=4\epsilon\mu(1-\mu)\mbox{,}
\end{align*}
the critical rate is $r_c(\epsilon)=1+2\epsilon$, and the connecting
orbit has the form $(x,\mu)=(-1+\mu(2-r/\epsilon),\mu)$. \pr{We will keep $\lambda_{\max} = 6$ fixed and vary $\rho$ between $0$ and $1/6$ which simultaneously varies $\epsilon$ between $0$ and $0.25$ and $r$ between $0$ and $1.5$. Therefore we will always have $\epsilon\ll 1$ corresponding to gentle but long ramps.}

% We will consider adding white noise to the dynamics of $x$ \eqref{sn ODE}, with the aim to create noise and rate-induced tipping according to a certain probability. This probability will depend on the rate of the shift of the system $\epsilon$, and the intensity of the noise.

The dynamics of $x$ for the system \eqref{sn ODE}--\eqref{lambda ODE},
modified by adding noise to \eqref{sn ODE}, are described by a
stochastic differential equation (an example of the general
equation~\eqref{gen SDE}):
\begin{equation}
\mathrm{d}X_{t} = [(X_{t} + \lambda(t))^{2} - 1]\mathrm{d}t + \sqrt{2D}\mathrm{d}W_{t}\mbox{.}
\label{Langevin}
\end{equation} 
The expressions (for example \eqref{eig prob}) for probability of
escape refer to a strip of half-width $\delta$ around a deterministic
reference trajectory $\tilde x(t)$. We choose the trajectory
$x^{u}(t)$ (blue curve in Figure~\ref{TP eps small} on the unstable
manifold $W^{u}(S_{-})$) and a strip with a fixed width $2\delta = 3$
around $\tilde x(t)$. Then we use the approximations for \prnew{$\mathbb{P}$, \eqref{Modes prob}, \eqref{eig prob}, \eqref{Flux prob}}
derived in section~\ref{sec: Methods} to
find % .then subtracting this away from
% $1$ will give
the probability of escape.

As described in section~\ref{sec: Methods}, we use co-moving coordinates

\begin{equation*}
y(t) = x(t) - x^{u}(t)\mbox{,}
\end{equation*}

\noindent such that the domain for $y$, $[-\delta,\delta]$ is fixed for
all $t$. This transformation alters the ODE given in \eqref{sn ODE}:

\begin{align}
\label{coord change initial ode}
\dot{x} &= f(x,\lambda) = (x+\lambda)^{2} - 1\mbox{,} \\
\nonumber
\dot{y} + \dot{x}^{u} &= f(y+x^{u},\mu)\mbox{,}  \\
\nonumber
\dot{y} &= (y + x^{u} + \lambda)^{2} - 1 - \dot{x}^{u}\mbox{,} \\
\nonumber
&= y^{2} + \underbrace{2(x^{u} + \lambda)}_{\mbox{$-c_{1}(t)$}}y + \underbrace{(x^{u} + \lambda)^{2} - 1 - \dot{x}^{u}}_{\mbox{equals\ $0$}}\mbox{,}
\end{align}  

\noindent and so, instead of \eqref{sn ODE}, we can express the new ODE as:

\begin{equation}\label{eq:gexample}
\dot{y}(t) = yg(y(t),t) = y^{2}(t) - c_{1}(t)y(t)\mbox{,}
\end{equation} 

\noindent where $c_{1}(t)$ is a time dependent scalar. Thus,
$yg(y,\mu(t))$ is used in the eigenvalue problem of the
Fokker-Planck equation \eqref{eigenvalue FP}, with
Dirichlet boundary conditions:

\begin{equation*}
P(-\delta,t) = P(\delta,t) = 0\mbox{.}
\end{equation*}

\begin{figure}[h!]
        \centering
        \subcaptionbox{\label{Eigenvalue spectrum all t}}[0.45\linewidth]
                {\includegraphics[scale = 0.3]{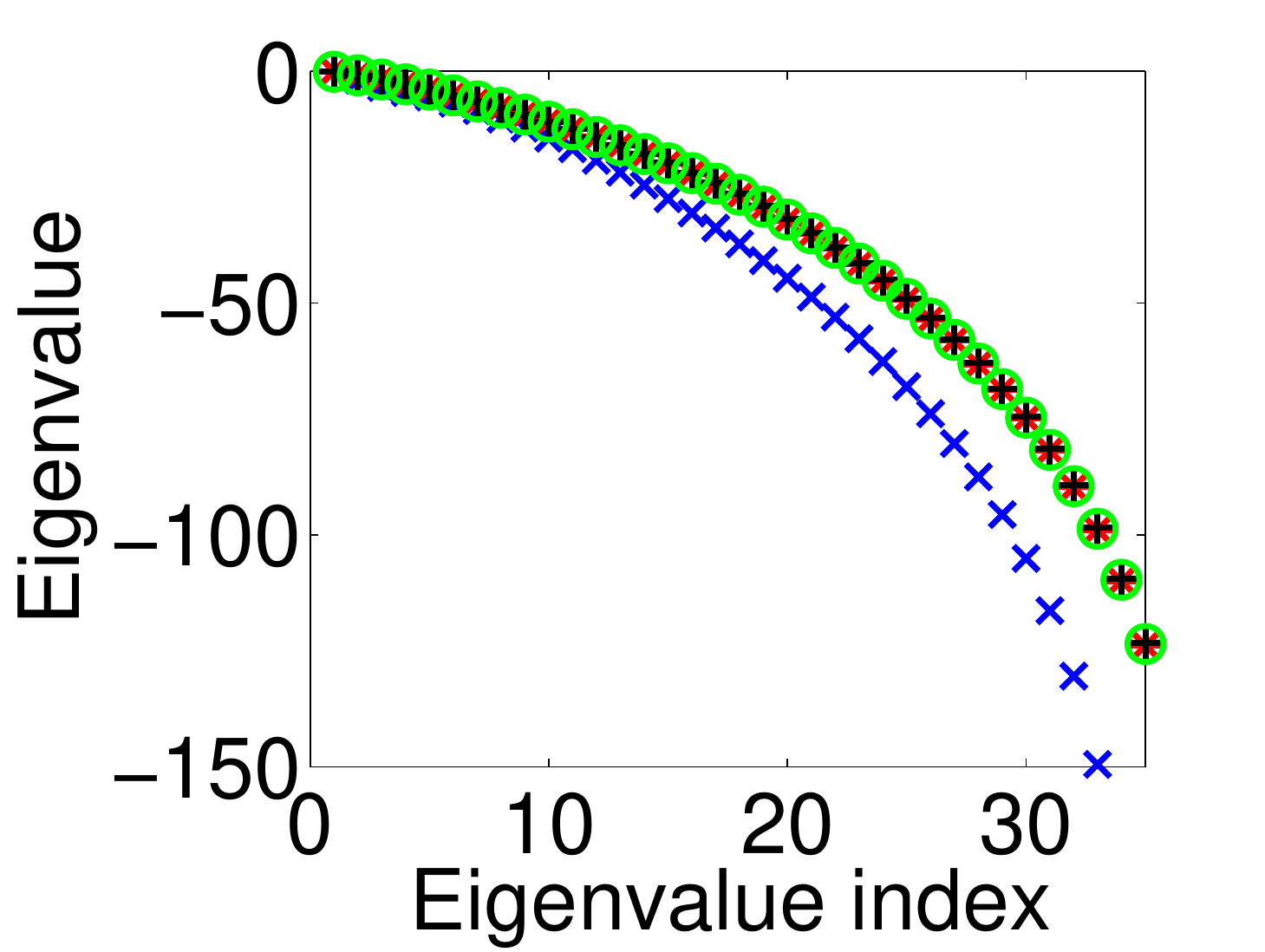}}
        \hfill %add desired spacing between images, e. g. ~, \quad, \qquad, \hfill etc.
        \subcaptionbox{\label{dominant eigenvalues all t}}[0.45\linewidth]
                {\includegraphics[scale = 0.3]{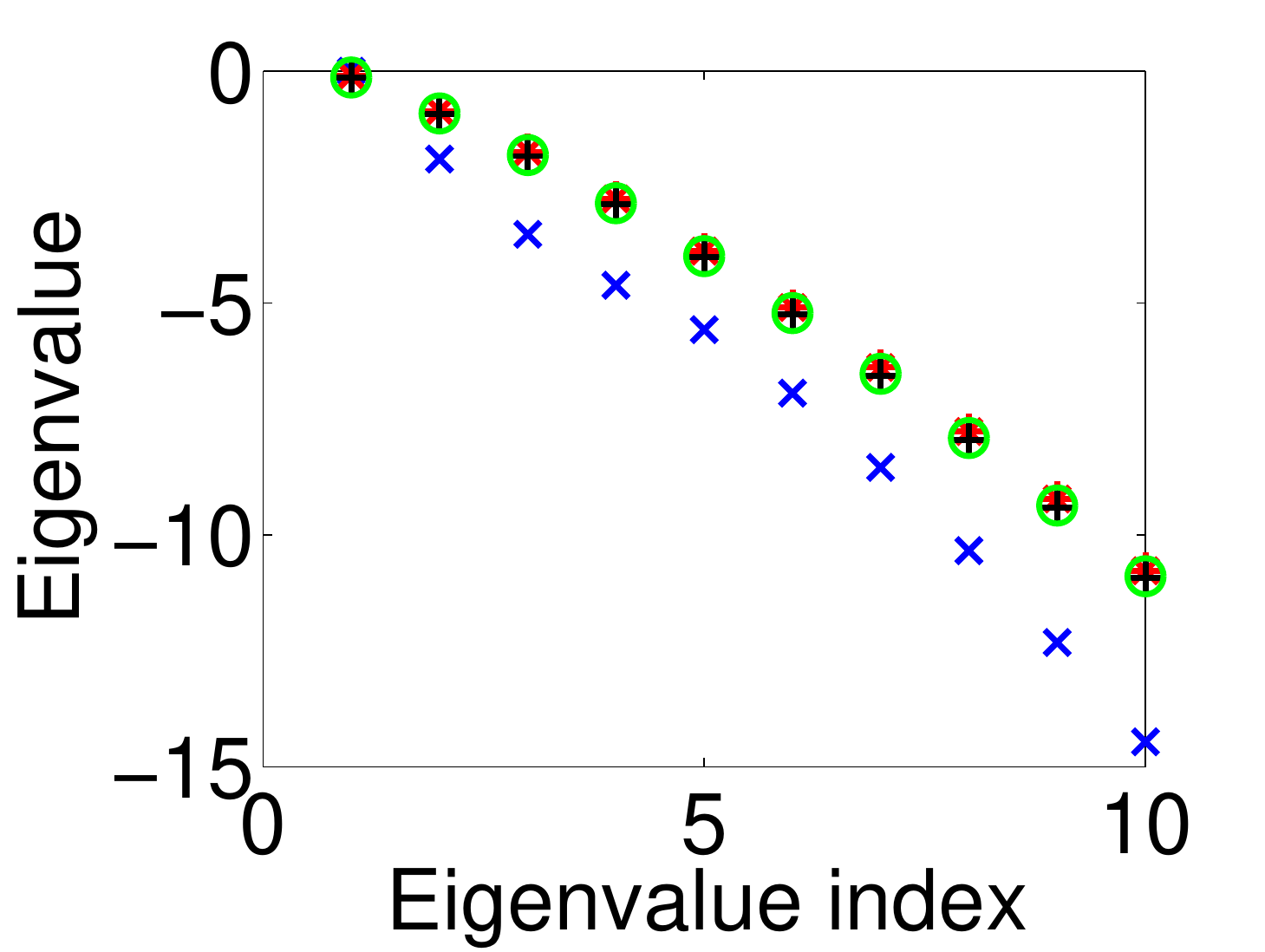}}
        ~ %add desired spacing between images, e. g. ~, \quad, \qquad, \hfill etc.
        \caption{Eigenvalue spectrum of the linear operator $A(t)$ for the Fokker-Planck equation corresponding to the system \eqref{sn ODE}--\eqref{lambda ODE}, with \pr{$\rho = 0.14$ ($\epsilon=0.21$, $r=1.26$),} $D = 0.06$ and at times $t = -10$ (blue cross), $t = 0$ (red star), $t = 1$ (black plus) and $t = 2$ (green circle).}\label{Eigenvalues all t}
\end{figure}

\noindent Figure~\ref{Eigenvalue spectrum all t} gives the spectrum of
the eigenvalues for different fixed times $t$ and for noise level $D =
0.06$. The eigenvalues for the times when the system \eqref{sn
  ODE}--\eqref{lambda ODE} is close to stationary (roughly all $t
\notin [-3,3]$), is given in
blue. %This is the same set of eigenvalues as presented in Figure \ref{Eigenvalues}.
The eigenvalues during the ramp (when $\dot\lambda$ is of order $1$)
are given by the other colors, namely at $t = 0$ (red star), $t = 1$
(black plus) and $t = 2$ (green circle).
%\js{\textbf{Nearly invisible in the
%    picture, may-be, use different symbols for each set. Eg, o, x and
%    +}} 
All sets of eigenvalues
$\gamma_k$ are on parabolas ($\gamma_k=\re \gamma_k\sim -k^2$)
(Figure~\ref{Eigenvalue spectrum all t}), with a nearly linear
relationship $\gamma_k\sim-k$ for the dominant eigenvalues (Figure~\ref{dominant eigenvalues all t}) and $\gamma_{1}(t) \approx 0$ for
all $t$.

Notice, that the first three eigenvalues during the ramp 
%\js{I changed
%  away from ``quasi-stationary'' here and above as it is confusing.}
are all greater than the second eigenvalue for times when the system
is close to stationary, for example at $t = -10$ (blue). Thus, the
contribution of additional modes is more significant during the period
of the
shift. % To achieve the same level of accuracy for a quasi-stationary density (for \pr{$\rho = 0.14$}) compared with a single mode approximation for a close to stationary density, three modes is required.}
The computational study in Section~\ref{sec: Probability} will compare
single- and three-mode approximations. % , which should provide a more
% accurate approximation to the probability density and probability of
% escape.

We return to Figure~\ref{Density comparisons} which displays different
types of densities to be compared with the \prnew{reference} probability density
from simulations in blue. We will focus on the comparison between the
single- (green (light gray)) and three-mode approximations (red (dark gray)) to the \prnew{reference} probability
density \pr{for ramping speed $\rho = 0.14$ ($\epsilon = 0.21$, $r =
  1.26$)}. According to our analysis all differences between
approximations and simulation results are caused by the non-zero
$\epsilon$. % The left column corresponds to a small ramping speed $\rho =
%0.2$ ($\epsilon=0.15$, $r=0.45$) and the right a moderate ramping
%speed, $\rho = 0.7$ ($\epsilon=0.525$, $r=1.62$). 

%For a slow ramping speed, both the single mode and first three modes give an accurate representation of the true probability density for all times $t$ (Figure \ref{Density comparisons} left column). Therefore, we will expect a good match for the probability of escape calculated using either a single mode or three modes with the true escape for slow ramp shifts $\rho$. Let us now consider a moderate ramping speed, $\rho = 0.7$ ($\epsilon=0.525$, $r=1.62$) and examine how well the modes approximate the true density, see panels (d), (f) and (h). 

Initially both the single-mode and three-mode approximation match the probability density well (not shown), however, a visible deviation appears in the single-mode at $t = 0$, see panel \pr{(b)}. The single-mode approximation develops a larger tail than the density from simulations (and the three-mode approximation). \change{This corresponds to an overestimation of the escape and hence, the peak of the density is lower.} \pr{The density has also shifted further along the $x$-axis because the density instantaneously adjusts to the effective potential, as previously discussed.} In contrast the three-mode approximation is still providing a good match to the \prnew{reference} probability density. %As the densities shift further along the $x$-axis there is a brief period at $t = 1$ when the three densities all coincide with one another, panel (f).\pr{At $t = 1$ (panel (c)), %there is once again a difference between the single mode and true density. \change{However, 
For larger times the single-mode is underestimating the escape compared to the simulations. \pr{The single-mode approximation converges back to the \prnew{reference} density for $t>2$ (panel (d)) when the shift slows down.} % These results would tend to suggest that the single mode is more responsive to the shift than the simulations density.
% and there exists a crossover at $t = 1$ when the densities coincide (the simulations density is widening and the single mode is narrowing). 
The three-mode approximation \pr{% nearly all the way through
  follows the \prnew{reference} density with greater accuracy
  throughout.} %apart from in panel (h). This corresponds to a time $t = 2$, where a small error develops in the tail of the 3 modes, an overestimation of escape. For greater accuracy one would need to use more modes. \change{All modes shortly after $t = 2$ converge and remain centred at $x = -4$ for the remaining time span, with slight differences in the peak heights reflecting different escape probabilities.}

% We continue with a systematic comparison of the probabilities for
% escape given by the single-mode approximation \eqref{eig alysis on the accuracy of the single- and three-mode
% approximations for the probability of escape (corresponding to a noise
% and rate-induced tipping event). Furthermore, we evaluate the
% probability of escape calculated using the formulas we have for the
% leading eigenvalue and probability flux.

\section{Systematic parameter study in $\rho$ and $D$}
\label{sec: Probability}

In this section, we will compare different approximations for the
probability of escape (noise and rate-induced tipping) in the two
parameter $(\rho,D)$ - plane. We will use
Monte-Carlo simulations (described below) as the reference for the
probability of escape.  Section \ref{sec: Methods} proposed two
approximate expressions for the single-mode approximation:
\begin{enumerate}
\item[(a)] the escape probability based on %the quasi-stationary density $P_*$
  %(given in \eqref{eq:pstar}), which is 
  the probability flux $J =
  D/p_{12}(a)$ with $a=-\delta$\prnew{, given in \eqref{Flux prob}};
\item[(b)] the escape probability based on the first-order approximation of
  the leading eigenvalue $\gamma_1$, given in \eqref{eig prob} (also
  with $a=-\delta$).
\end{enumerate}
We also compute the escape probabilities based on single- and
three-mode approximations, by solving the ODE eigenvalue
problem \eqref{eigenvalue FP} with Dirichlet boundary conditions at
$-\delta$ and $\delta$ numerically for the first \prnew{$n$} modes
\prnew{($(\gamma_k,v_k)$, $k=1,\ldots,n$)} together with the ODE
\eqref{akdot}, such that the escape probability \prnew{is given by \eqref{Modes prob} for} %equals$1-\int_{-\delta}^\delta P_k(x,T_\mathrm{end})\d x$ (for 
$n=1$ and
$n=3$. % with which we will compare to the approximations using the
% first three modes and first mode of the linear operator of the
% Fokker-Planck equation. However, calculating the modes numerically is
% equivalent to solving the Fokker-Planck equation a PDE and so we will
% also use the two formulas derived earlier in \change{Section \ref{sec:
%     Methods}}. Namely the probability flux $J = D/p_{12}(a)$ and
% leading eigenvalue $\gamma_{1}(t)$, equation \eqref{lead eig}, which
% give approximations for the leading single mode and hence, also give
% an approximation to the true probability of escape.

\begin{figure}[h!]
        \centering
        \subcaptionbox{\label{Sim prob}}[0.45\linewidth]
                {\includegraphics[scale = 0.3]{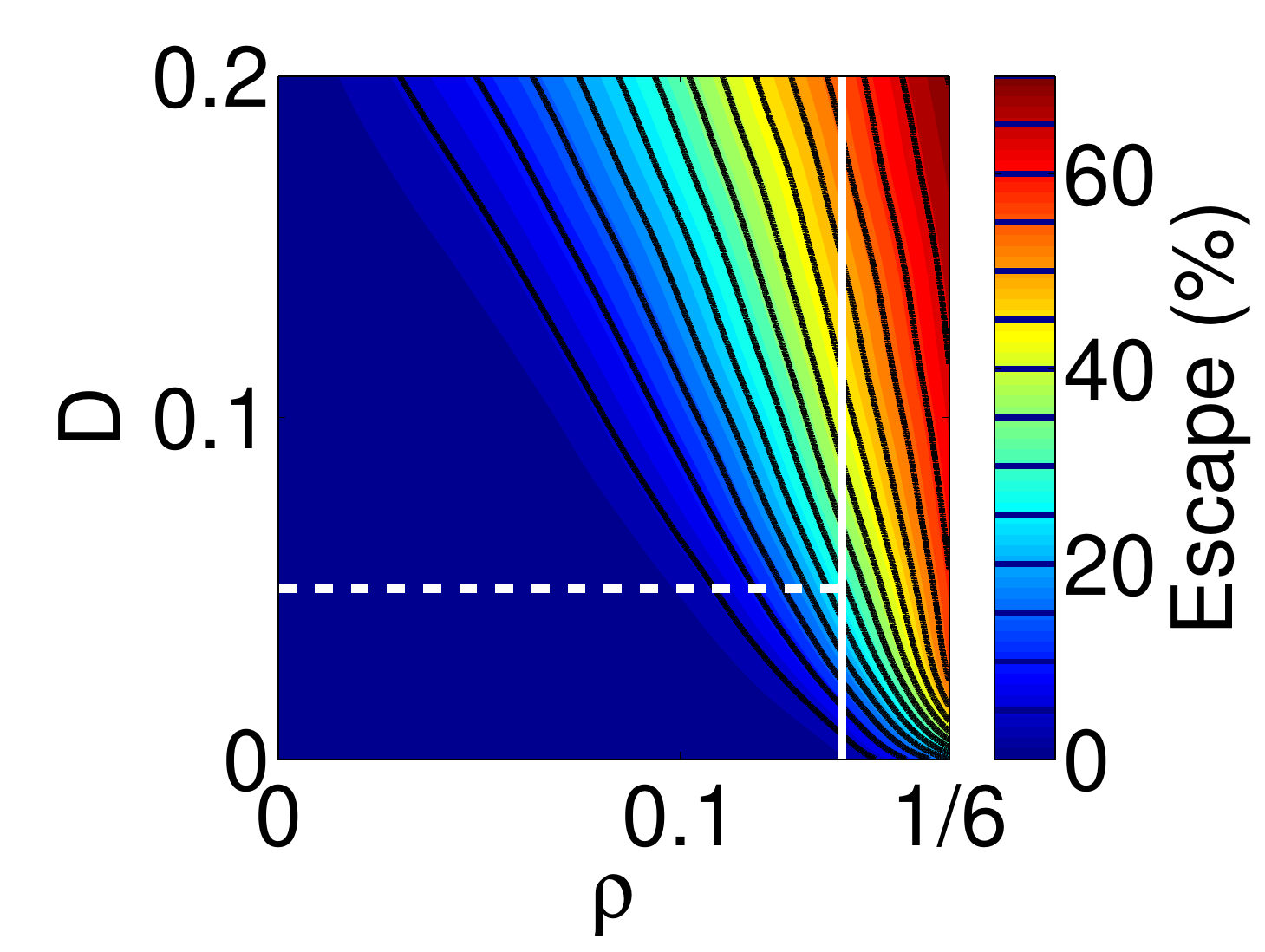}}
        \hfill %add desired spacing between images, e. g. ~, \quad, \qquad, \hfill etc.
        \subcaptionbox{\label{c1 plot}}[0.45\linewidth]
                {\includegraphics[scale = 0.3]{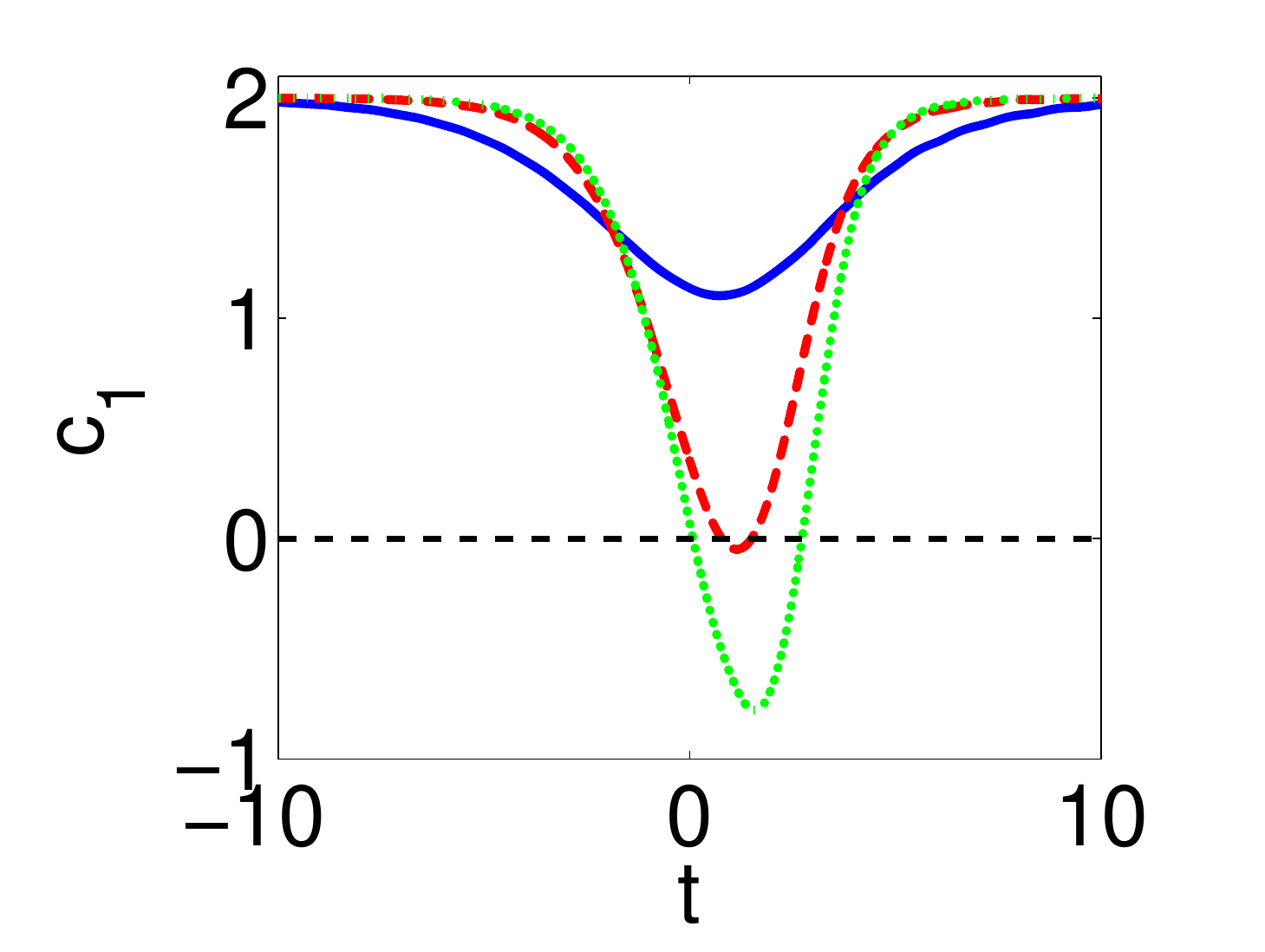}}
        ~ %add desired spacing between images, e. g. ~, \quad, \qquad, \hfill etc.
        \caption{(a) Overview probability of escape in the $(\rho,D)$
          - plane using simulations (with smoothing), where contours
          are spaced for every $5\%$ of escape. Vertical white line
          indicates value of $\rho$ such that $\min(c_{1}) = 0$ and
          horizontal dashed line shows lowest value of $D$ for which
          the probability can be calculated for the modes. (b) Time
          profile of the critical points of $U(y,t)$, one of them
            is always $y=0$, the other is $c_{1}(t)$ for $\rho =
          0.08$ (blue solid), $0.14$ (red dashed) and $0.16$
          (green dotted).}\label{Overall prob}
\end{figure}

Figure~\ref{Sim prob} shows the probability of noise and rate-induced
tipping occurring in the two parameter $(\rho, D)$ - plane, calculated
using Monte-Carlo simulations. This has been performed by starting
with a large number of realizations at $x_{0} = -1$ at $t_{0} = -10$
and evolving according to the SDE \eqref{gen SDE}. The fraction of
realizations that pass $x_{T} = 4$, and, hence, go to $+\infty$ in
finite time, approximates the probability of tipping (or probability
of escape). The reference
probability is not derived from the number of realizations escaping
the strip $\{y(t)\in[-\delta,\delta]: t\in[t_0,T_\mathrm{end}]\}$ but
by the fraction of realizations that have crossed an arbitrary line
$x_{T} = 4$ (the choice of $\delta$ and $x_T$ is such that this
difference has a small effect).

Figure~\ref{Sim prob} shows the probability of escape (in $\%$) for
all ramping speeds $\rho$ up to $\rho_{c} = 1/6$ and a range of noise
levels $D$. The color contours indicate that the probability of escape
is small for small $\rho$ and $D$. As $\rho$ increases towards
$\rho_{c}$ and the noise level increases so does the probability of
escape, reaching approximately $70\%$ probability of escape for $\rho
= \rho_{c}$ and $D = 0.2$.

\paragraph{Region in $(\rho,D)$-plane considered} We can expect
  the single-mode (or three-mode) approximations to be accurate
  only in a range of parameters $\rho$ up to a value
  $\rho_{\max}=0.14$ that is slightly smaller than the critical value
  $\rho_c=1/6$ (where tipping occurs without noise). The reason for
  this is in the error terms when replacing the dynamic Fokker-Planck
  equation for the density with its projection onto leading
  time-dependent mode(s). These error terms are only small if the time
  derivative of the (time-dependent) drift $yg(y,t)$ in
  \eqref{eq:shiftFP} is small.

As introduced in
\eqref{eq:shiftSDE} in Section~\ref{sec:approx}, the eigenvalue
problem for the Fokker-Planck equation is solved in a co-moving
coordinate system along the path $\tilde x(t)$ ($y=x-\tilde x(t)$)
such that the path is centered at $y = 0$ within the fixed domain
$y\in[-\delta,\delta]$ and the drift is given by
\begin{math}
yg(y,t)=y(y - c_{1}(t))
\end{math}
(where $c_1(t)$ is given in \eqref{coord change initial ode}).  Its
potential $U(y,t) = -\int y(y-c_1(t))\mathrm{d}y$ has a well at $y =
0$ and a hill top at $c_1$ for $c_{1}> 0$, but a hill top at $y=0$
(and a well at $c_1$) for $c_{1} < 0$.
% We choose the range of For negative
% $c_{1}$
% we may incur escape to the left boundary if our domain (that is,
% $\delta$) is too small. Therefore a natural cut off would be to
% determine a $\rho$ for which the minimum value of $c_{1}$ is zero.
Figure~\ref{c1 plot} illustrates the time profile of $c_{1}(t)$ for
different values of the drift speed $\rho$. The limit of $c_1(t)$ for
$t\to\pm\infty$ is $2$ such that $c_1(t)\approx2$ for $t$ close to
$t_0=-10$ and $T_\mathrm{end}=10$. For small drift speeds $\rho$ the
deviations of $c_1$ from its asymptotic value are small (blue solid curve in
Figure~\ref{c1 plot}), while for $\rho = 0.16$ $c_1(t)$ becomes
negative for some time interval, making the trajectory $x(t)=\tilde
x(t)$ or $y(t)=0$ locally repelling (green dotted curve in Figure~\ref{c1
  plot}). 

The error of the single-mode approximation is small if $|\dot c_1(t)|$
is small, which is the case for $t$ near $t_0$ and $T_\mathrm{end}$,
and for $t$ close to the minimum of $c_1$. If $c_1(t)>0$ for all $t$
then the time $t_{\min}$, when $c_1(t)$ is minimal, correspond to
those times where escape is most likely to occur, since at these times
the potential barrier is smallest. At times near $t_{\min}$ the mode
approximation error is also small since $|\dot c_1(t)|$ is
small. 

However, if $c_1(t)<0$ for a range of $t$ then escape occurs with a
non-small probability at times when $|\dot c_1(t)|$ is not small, \prnew{leading to an error in the single-mode approximation that is not small.}

% the system shifts relatively slowly
% and this is reflected by the slow change and small range of values of
% $c_{1}$. For a faster shift, \pr{$\rho = 0.16$} (green) the system
% remains stationary for longer ($c_{1} = 2$) but shifts quicker
% corresponding to a fast change over a wider range of $c_{1}$
% values. Notice that i
Hence, we choose a range for the parameter $\rho$ such that $c_1(t)$
stays positive along the entire path for all $\rho$.

% Even though $c_{1}$ is negative for a short period of
% time (making the trajectory $x(t)=\tilde x(t)$ or $y(t)=0$ locally
% repelling) tipping does not occur without noise as $\rho$ is still
% smaller than the critical value $\rho_c=1/6$. We find that \pr{$\rho =
%   0.14$} (red) corresponds to the minimum of $c_{1}$ being roughly
% zero and so we will use this as our cut off value as indicated by the
% vertical white line in Figure~\ref{Sim prob}. \pr{Notice that we can
%   consider cases of close encounters of rate-induced tipping by having
%   a long (large $\lambda_{\max}$) but slow (small $\rho$) shift.}

We also remove small values of $D$ ($D<0.05$) from our consideration,
since for $\rho<\rho_c$, but not close to $\rho_c$, escape
probabilities are small compared to errors in Monte-Carlo simulations
and in the numerical computations of the integrals needed for
$\gamma_1$ in \eqref{lead eig}. In this region the probability of
escape is exponentially small in $D$ (that is, of order $\exp(-C/D)$
for some constant $C>0$). % In principle we can consider everywhere
% between \pr{$0 \leq \rho \leq 0.14$}, though we also restrict our
% domain to above the horizontal white dashed line (Figure \ref{Sim
%   prob}) when calculating the modes. One reason is the simulations
% illustrate that the probability of escape is very small. The other
% reason is that it would be a lot more computationally expensive as we
% would require a finer grid spacing for small values of $D$.

\begin{figure}[h!]
        \centering
        \subcaptionbox{\label{True}}[0.45\linewidth]
                {\includegraphics[scale = 0.3]{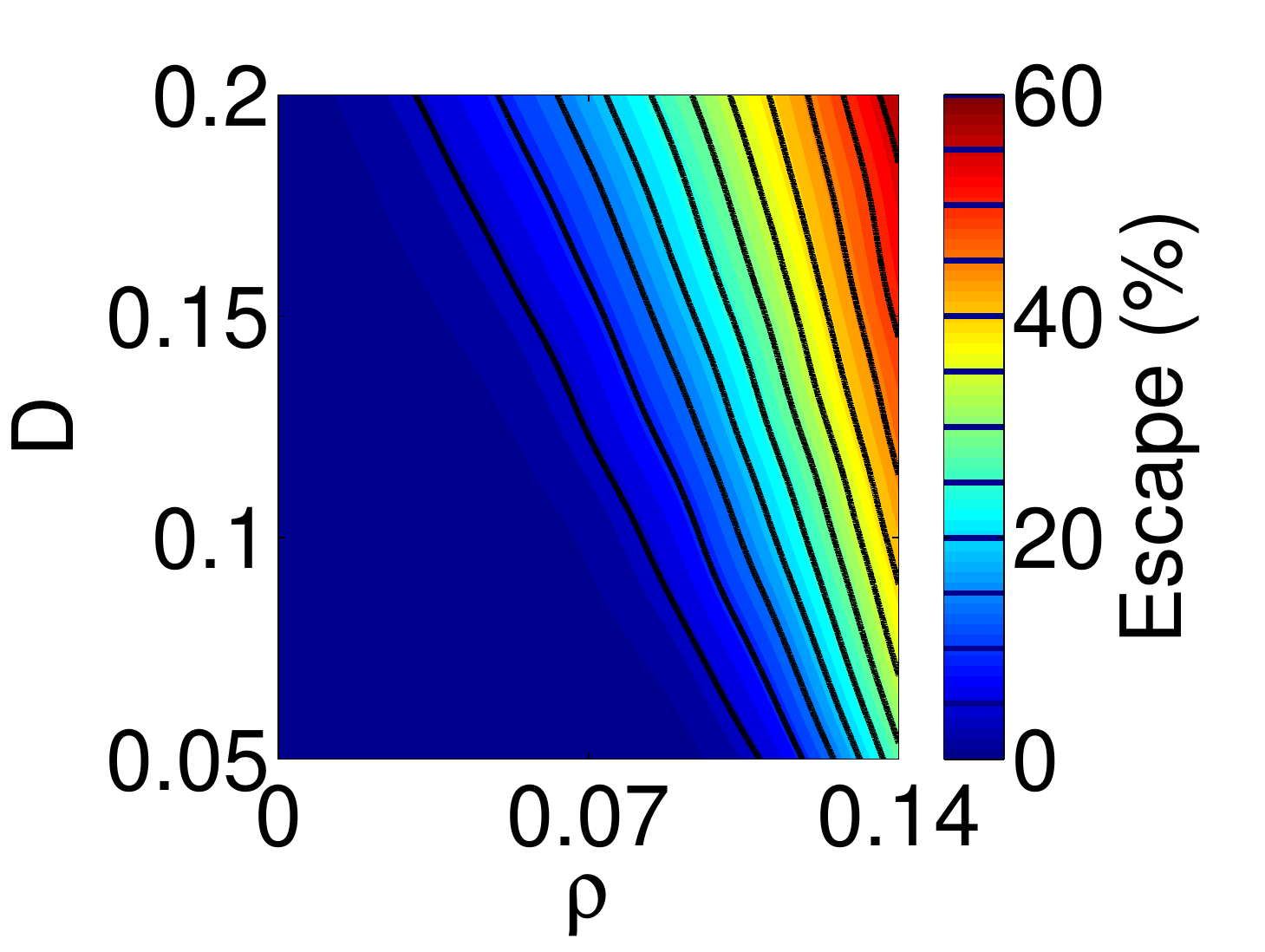}}
        \hfill
        \subcaptionbox{\label{3 modes}}[0.45\linewidth]
                {\includegraphics[scale = 0.3]{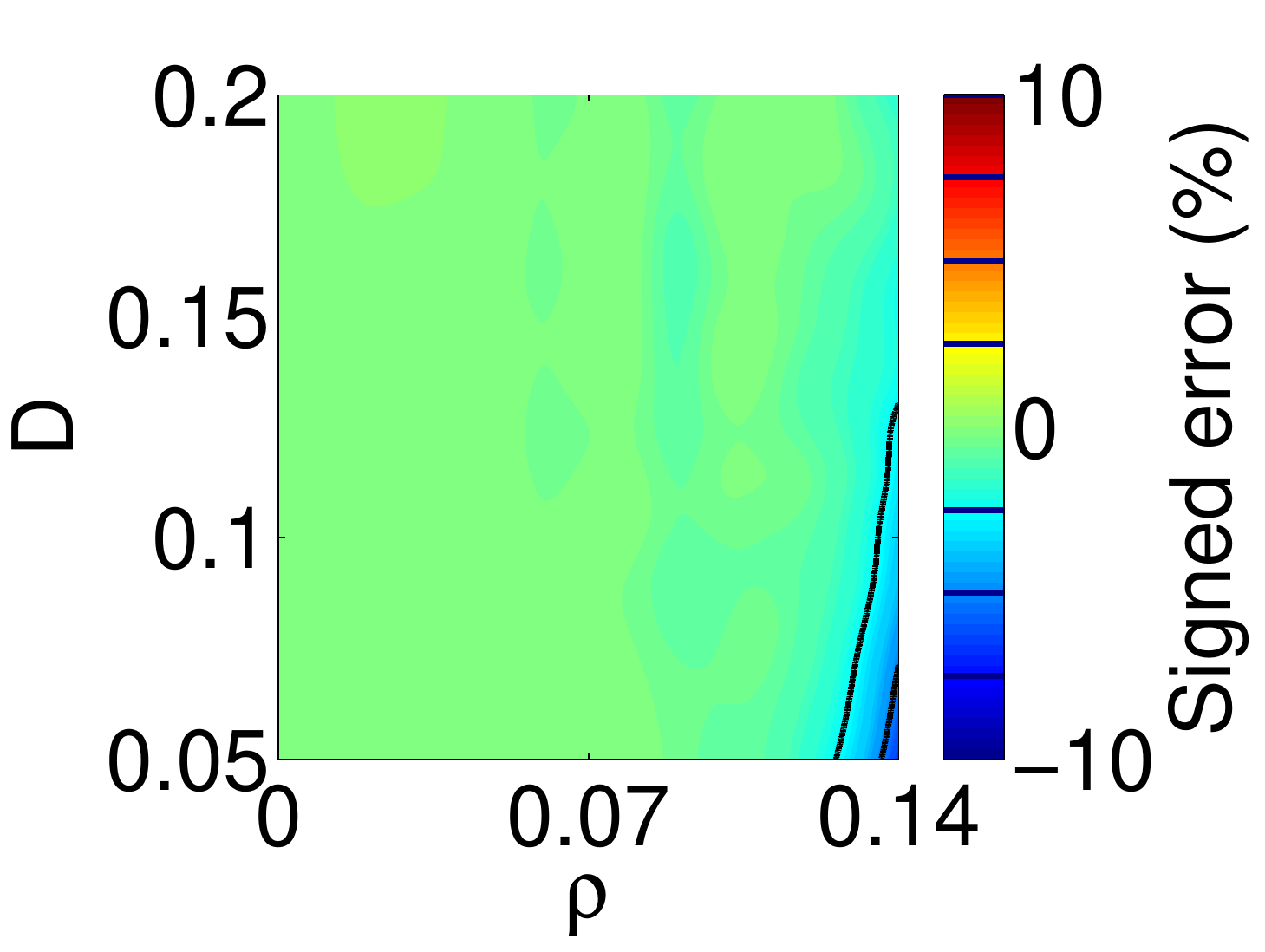}}
        ~ 
        \\
        \subcaptionbox{\label{formula 1 mode}}[0.45\linewidth]
                {\includegraphics[scale = 0.3]{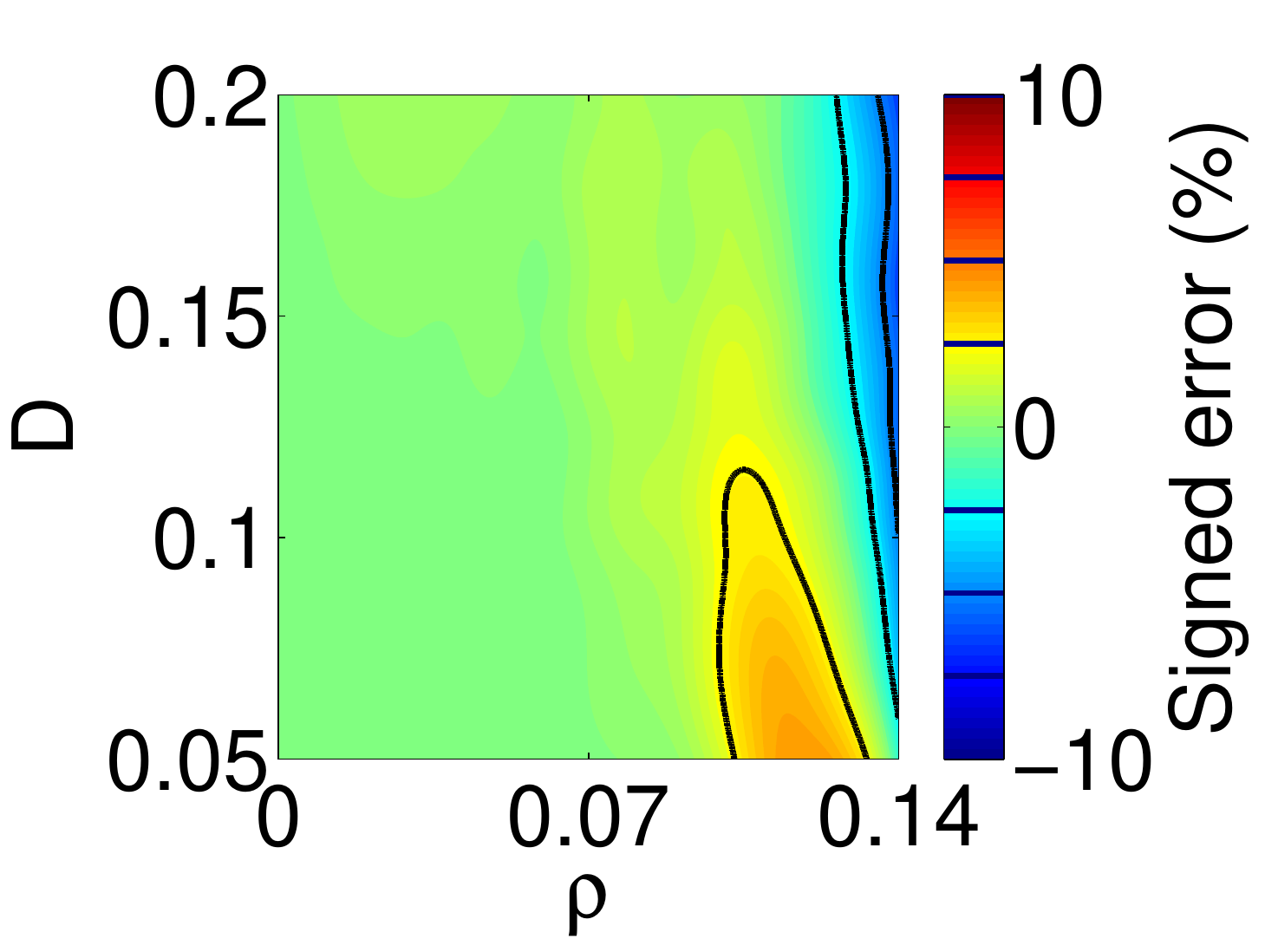}}
        \hfill
        \subcaptionbox{\label{kramers 1 mode}}[0.45\linewidth]
                {\includegraphics[scale = 0.3]{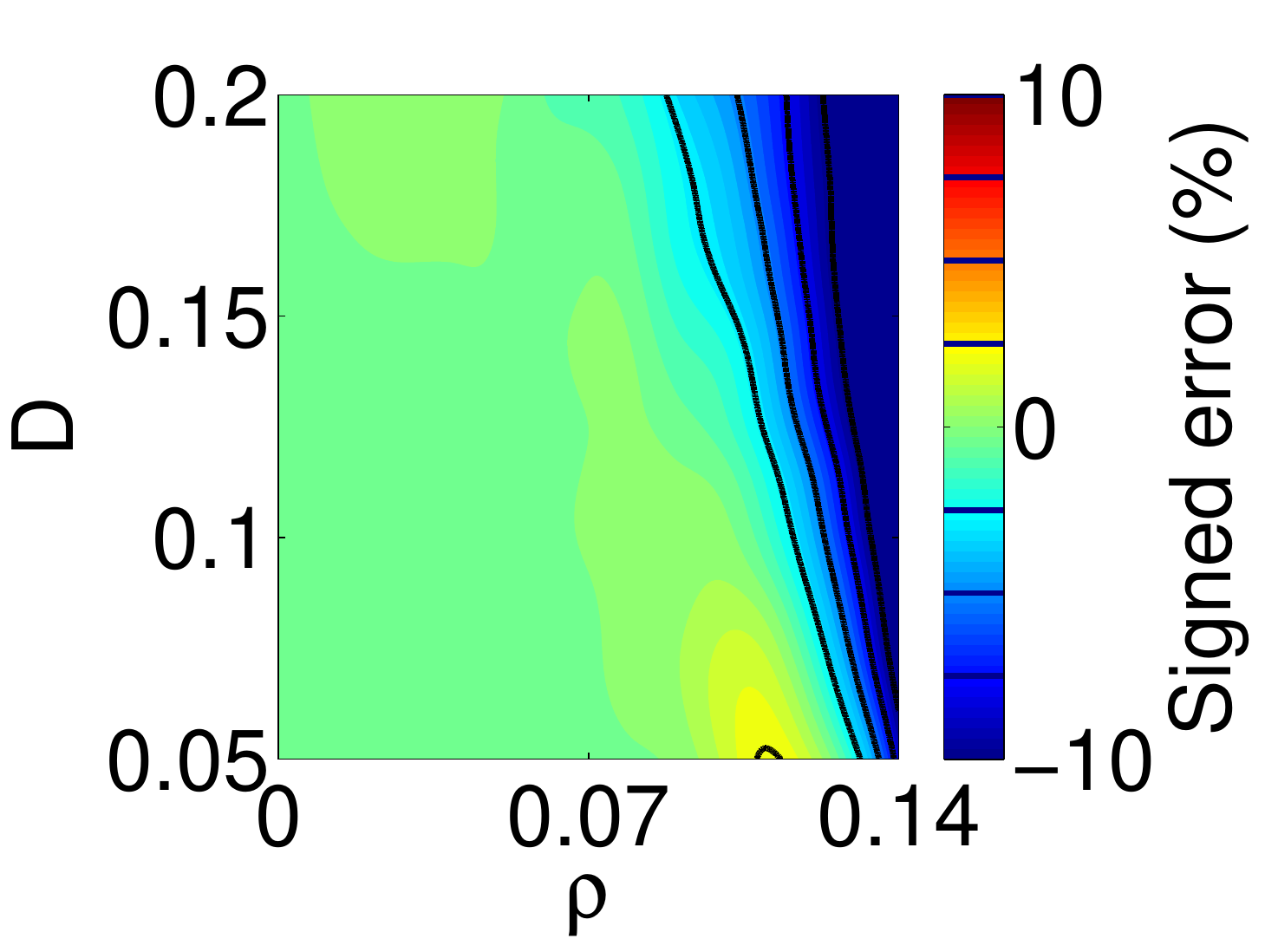}}
        \caption{(a) Reference probability of escape in $(\rho,D)$ - parameter plane, obtained using Monte-Carlo simulations (with smoothing) - observing the fraction of realizations that escape the potential landscape. (b)-(d) Evaluating approximation methods for probability of escape by plotting in color the signed error $\%$ between the approximation and the reference escape. Positive error (red) reflects an overestimation and negative (blue) an underestimation. Approximation methods used: (b) Numerically calculate first three ($n=3$) modes of the linear operator of the Fokker-Planck equation \eqref{Modes prob}. (c) Analytical single-mode approximation \eqref{eig prob} and (d) analytical probability flux approximation \eqref{Flux prob}. Contours are spaced at $5\%$ intervals for panel (a) and $2.5\%$ intervals for panels (b)-(d) with the zero contour omitted.}\label{Mode approximations}
\end{figure}

Figure~\ref{True} \prnew{shows the reference probability - the probability
of escape calculated using Monte-Carlo simulations for this restricted region} (with a slightly different
color scale to Figure~\ref{Sim prob}). The remaining three panels of
Figure~\ref{Mode approximations} give the signed error of the
approximation, compared to the reference, in percent. In the color
scale for these panels, a % measurWe would like
% to highlight that the color scheme is different to \ref{True} and so
% a
green (light gray) color represents good agreement between the approximation and
reference escape. A positive error (red) means an overestimation and a
negative error (blue) corresponds to an underestimation when using the
approximation method.

%\begin{figure}[h!]
%        \centering
%        \subcaptionbox{Single mode\label{1 mode}}[0.45\linewidth]
%                {\includegraphics[scale = 0.3]{rtip_escape_grid_1mode.eps}}
%        \hfill %add desired spacing between images, e. g. ~, \quad, \qquad, \hfill etc.
%          %(or a blank line to force the subfigure onto a new line)
%        \subcaptionbox{3 modes\label{3 modes}}[0.45\linewidth]
%                {\includegraphics[scale = 0.3]{rtip_escape_grid_3modes.eps}}
%        ~ %add desired spacing between images, e. g. ~, \quad, \qquad, \hfill etc.
%          %(or a blank line to force the subfigure onto a new line)
%        \caption{Comparison of using a single mode or 3 modes to approximate the probability of escape in the 2 parameter $(\epsilon,D)$ - plane}\label{Mode approximations}
%\end{figure} 

Figure~\ref{3 modes} shows the probability of escape calculated using \eqref{Modes prob} with the first $n=3$ instantaneous eigenmodes of the linear operator of the Fokker-Planck equation. % This panel providesis designed to highlight the potential of the method (probability of escaping a strip $S_{\delta}$ surrounding a path $\tilde{x}(t)$) and will offer as a reference for the analytical approximations for this method.
  The three-mode approximation was computed by solving the ODE eigenvalue problem \eqref{eq:FPeig} numerically. It approximates the reference probability of escape over the specified region well, except for $\rho$ close to $0.14$ and small noise levels. % does three-mode approximation give a small underestimation of the reference probability  for escape.
%, indicating further modes may be required. %Using just the first instantaneous eigenmode offers a very good agreement with the true escape for small \pr{$\rho<0.1$}, see Figure \ref{1 mode}. The single-mode then generally provides an overestimation with a larger error for smaller noise levels $D$.

Figures \ref{formula 1 mode} and \ref{kramers 1 mode} compare the single-mode approximation \eqref{eig prob} and the approximation using the probability flux \eqref{Flux prob} to the reference escape. The single-mode approximation offers a very good agreement with the reference escape for $\rho<0.1$. For larger $\rho$ the formula gives an overestimation of the escape for small noise levels. %The error from the formula for the probability flux has a similar structure to the perturbation formula error but surprisingly gives a better approximation of the true escape on the whole than the single mode. 
The probability flux again approximates the probability for small $\rho$ values well, \pr{but when the probability of escape increases to above $20\%$ the probability flux underestimates the reference escape.} % and even though there is an overestimation for larger $\rho$ and small $D$ it is not as high as that seen for the formula and/or numerically calculated single mode. It is only when the escape is largest (large $\rho$ and $D$) does the perturbation formula provide a better approximation. However, the probability flux providing a better approximation than the perturbation formula is likely to be specific to this example and not true in general. 
% Both formulas are an approximation to the probability calculated numerically from the first eigenmode. Therefore we will develop this further by analyzing the error between the probabilities calculated numerically from the single-mode with the two formulas, see Figure~\ref{formulas 1 mode}.       

\begin{figure}[h!]
        \centering
        \subcaptionbox{\label{formula 1 mode error}}[0.45\linewidth]
                {\includegraphics[scale = 0.3]{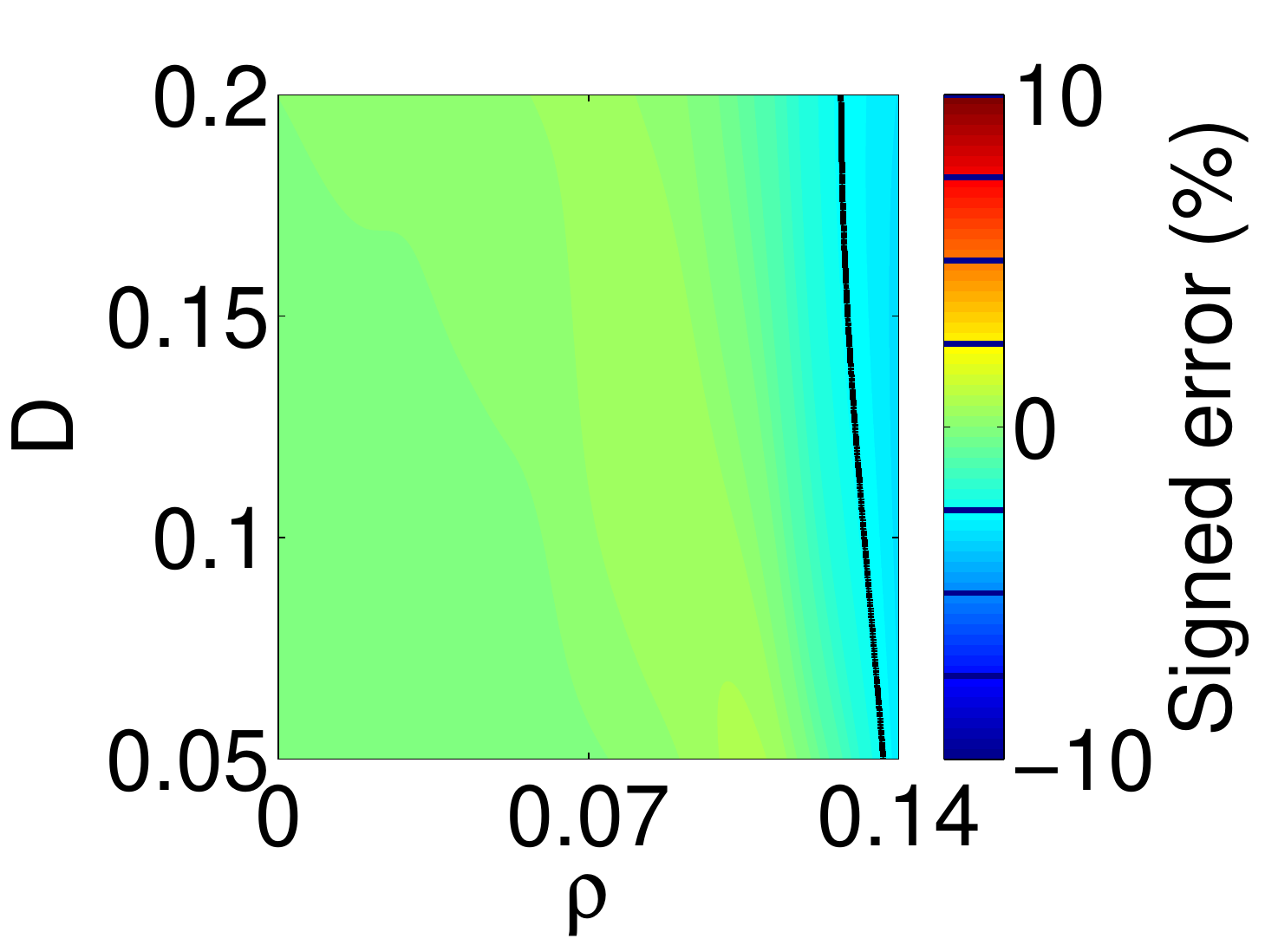}}
        \hfill
        \subcaptionbox{\label{kramers 1 mode error}}[0.45\linewidth]
                {\includegraphics[scale = 0.3]{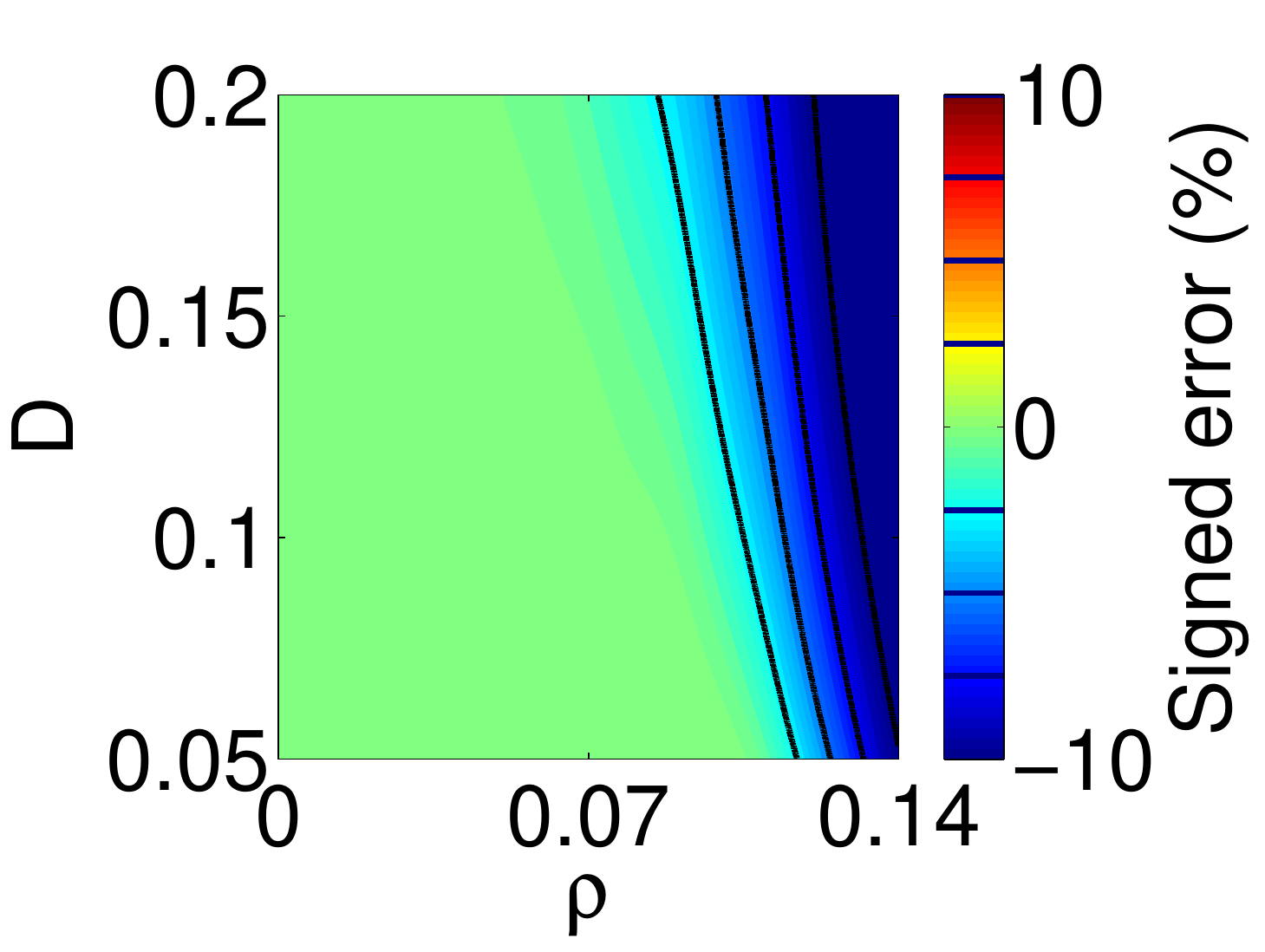}}
        \caption{Color plots of the signed error between the probability calculated numerically from the ($n=1$) single-mode \eqref{Modes prob} and (a) the single-mode approximation \eqref{eig prob} or (b) the probability flux approximation \eqref{Flux prob}. A positive error (red) corresponds to the prescribed approximation overestimating the numerical probability calculated from the single-mode, whereas a negative error (blue) represents an underestimation. Contours are spaced at $2.5\%$ intervals with the zero contour omitted.}\label{formulas 1 mode}
\end{figure} 

Figure~\ref{formula 1 mode error} \prnew{compares the single-mode approximation formula \eqref{eig prob} to the exact single-mode approximation, solving \eqref{eq:FPeig} and \eqref{Modes prob} for $n=1$ numerically. We see that the difference
%between the numerical solution for the single-mode approximation and
%the linear perturbation formula 
is much smaller than
the error caused by the approximation of the Fokker-Planck equation
with a single mode, see Figure~\ref{formula 1 mode}.} % gives an
% extremely good match to the single mode over the two parameter region
% considered. The perturbation formula only slightly underestimates the
% probability calculated from the single mode for $\rho$ close to
% \pr{$0.14$}.

In contrast, the probability flux \eqref{Flux prob} systematically
underestimates the escape probability when it is greater than
\pr{$20\%$}; see Figure~\ref{kramers 1 mode error}. This is as
expected because the estimate \eqref{Flux prob} based on a spatially
  constant probability flux assumes the flux escaping at the
  boundary $+\delta$ re-enters at the boundary $-\delta$. % realisations
% that escape are put back in at the left boundary. This is in contrast
% to substituting them back in according to the distribution of the
% remaining realisations.
The difference between the two estimates is larger
\pr{when the escape probability is high}, because the assumption of spatially constant flux is only approximately true if escape is sufficiently rare compared to the time it takes for realizations to reach the potential well from the boundary at $-\delta$. % put back in
% according to the distribution are far more likely to escape earlier
% than those put in at the opposite boundary.
This leads to an underestimate for those values of $D$ and $\rho$ when
the probability of escape is greatest.

\section{Discussion}
\label{sec: Discussion}

We have provided approximations for the critical rate for
deterministic rate-induced tipping and for the probability of
noise-induced tipping during parameter shifts. These
  approximations are valid for parameter shifts that are asymptotic to
  constant parameter values for $t\to\pm\infty$ and that are ``long
  but gentle'': the small parameter $\epsilon$ is the ratio between
  maximal ramp (shift) speed $r$ and the length of the parameter
  shift. The deterministic critical rate $r_c(\epsilon)=r_0+\epsilon
  r_1+O(\epsilon^2)$ is an order-$\epsilon$ perturbation from the
  critical rate $r_0$ for constant-speed parameter shifts (which were
  discussed in
  \cite{ashwin2012tipping}). % but at slow speeds $\rho$ we have a
  % perturbation formula that approximates the probability of escape for
  % close encounters of rate-induced tipping
  % well.}
% \change{We also note that there is possibilities to approximate higher probabilities of tipping by considering examples of longer shifts over slower speeds.}
The approximation for the tipping probability in the presence of noise
is based on the instantaneous eigenmode expansion of the linear
operator for the Fokker-Planck equation to approximate the
quasi-stationary probability densities % . We are then able to calculate
% the probability of escape by evaluating the leamount of density that
% escapes
in a strip $S_\delta$ of half-width $\delta$ around the deterministic
trajectory (a connecting orbit between equilibria). Moreover, we have
derived a general perturbation formula to calculate the leading
eigenvalue, which approximates the probability calculated from a
single eigenvalue (the single-mode approximation) and thus, gives a
good approximation to the reference probability of tipping for small $\epsilon$.

% This paper complements the work of \citet{ritchie2015early}, which
% derived a boundary value problem for calculating the optimal path that
% can be used to help determine the most likely time of a noise and
% rate-induced tipping event. A natural progression after determining
% the most likely time of tipping is to investigate the probability of
% tipping. However, the probability for an event happening cannot be
% determined by the most likely path for escape as not all
% \pr{realizations} that escape will remain within a band width that
% distinguishes between trajectories that escape and those that do
% not. Instead we calculate the probability of following the
% deterministic trajectory (there is no tipping without noise) within a
% predetermined band width and subtract from $1$ to give the probability
% of escape. This requires knowing the quasi-stationary probability
% densities. However, to avoid solving the Fokker-Planck equation a PDE
% we use the instantaneous eigenmode expansion of the linear operator,
% satisfying a system of ODEs instead \eqref{akdot}.
The limitation of the proposed estimate using single-mode approximation is that it fails
  for \prnew{some maximal ramp speeds $r$ less than $r_c$, even for small noise levels $D$. A brief derivation in \review{Appendix \ref{app:crit rate expansion}} shows that the single-mode approximation is generally valid for parameters up to $r=r_0+\epsilon r_1/2$ such that the escape rate from the strip $S_\delta$ is maximal during times when the time-derivative of the shape of the underlying potential well is minimal. % See Supplementary Material at [URL will be inserted by publisher] for further discussion.
  }   

\begin{acknowledgments}
  P.D.L.R.’s research was supported by funding from the
EPSRC Grant No. EP/M008495/1, J.S. gratefully acknowl-
edges the financial support of the EPSRC via Grants No.
EP/N023544/1 and No. EP/N014391/1. J.S. has also received
funding from the European Union’s Horizon 2020 research and
innovation programme under Grant Agreement No. 643073.
\end{acknowledgments}

\begin{appendix}
\section{Expansion of the critical rate in the small parameter}
\label{app:crit rate expansion}

\noindent This section presents in more detail the expansion of the
critical rate for systems with a ramped parameter shift. The critical
rate is defined as the threshold at which a system fails to track the
continuously changing quasi-steady state and thus generating
rate-induced tipping.

The simplest example of rate-induced tipping is a system subjected to a linear parameter shift $r_\lin$, $\dot x=f(x+b\lambda)$ with $\lambda=r_{\lin}t$. We assume that the linearly shifted system ($y=x+b\lambda$)
\begin{equation}
\dot{y} = f(y) + r_{\lin}b
\label{lin system}
\end{equation}
($f:\R^n\mapsto\R^n$, $y(t)\in\R^n$) has a generic saddle-node
bifurcation at $r_{\lin} = r_0>0$, $y = y_0$. We assume that the
number of unstable dimensions of the equilibria changes from $0$ to
$1$ at the saddle-node and that the equilibria exist for
$r_{\lin}<r_0$. In this section we derive the first-order expansion of
the critical rate $r_c$ in $\epsilon$ for the system
\begin{align}
\dot{y} &= f(y) + b(r\Gamma(\mu)-r_0)\mbox{,} &
\dot \mu&=\epsilon \Gamma(\mu)\mbox{.}
\label{ydot}
\end{align}
%\noindent which sets $r_\lin = r\Gamma(\mu)-r_0$ in \eqref{lin system}. 
System \eqref{ydot} describes the scenario of a ramped shift $\dot
x=f(x+br\mu/\epsilon)$ with maximal speed $r$, again in shifted
coordinates $y=x+br\mu/\epsilon$. The direction of the shift is
determined by $b\in\R^n$ and $\Gamma(\mu)$ satisfies the following
properties
\begin{equation}\label{Gamma properties}
  \begin{aligned}
    &0=\Gamma(0)=\Gamma(1)\mbox{,}\\
    &1=\Gamma(\mu_\mathrm{crit})=\max\{\Gamma(\mu):\mu\in[0,1]\}, \,\mbox{(this defines $\mu_\mathrm{crit}$),}\\
        &0>\Gamma''(\mu_\mathrm{crit}), \,\mbox{(let $g_2:=-\Gamma''(\mu_\mathrm{crit})/2$),}\\
    &\Gamma(\mu)\in(0,1)\mbox{\ for all $\mu\in(0,\mu_\mathrm{crit})$
      and $\mu\in(\mu_\mathrm{crit},1)$}
  \end{aligned}
\end{equation}
such that $\Gamma$ has a unique non-degenerate maximum at
$\mu_\mathrm{crit}$. We make a change of coordinates $y = y_0 + v_0z$
to shift the origin in \eqref{ydot} to the saddle-node bifurcation at
$y = y_0$:
\begin{equation}
v_0\dot{z} = f(y_0 + v_0z) + b(r\Gamma(\mu)-r_0),
\label{ydot coord change}
\end{equation}
\noindent where $v_0$ is the right nullvector of $\partial f(y_0)$
(which has a one-dimensional nullspace by the assumption of a
saddle-node bifurcation at $y=y_0$ for \eqref{lin
  system}). Furthermore, this assumption implies
\begin{equation}
f(y_0 + v_0z) = \frac{1}{2}z^2\partial^2f(y_0)v_0^2 + \mathcal{O}(z)^3,
\label{Expansion simplified}
\end{equation}
Inserting \eqref{Expansion simplified} into \eqref{ydot
  coord change} and applying $w_0^T$ (the left nullvector of $\partial
f(y_0)$, scaled such that $w_0^Tv_0=1$) to \eqref{ydot coord change}
gives
\begin{equation*}
\dot{z} = a_0(r\Gamma(\mu)-r_0) + a_2z^2 + \mathcal{O}(z)^3,
\end{equation*}
where $a_0 = w_0^Tb$, $a_2 = \frac{1}{2}w_0^T\partial^2f(y_0)v_0^2$
are both non-zero by the assumption of a generic saddle-node. As we
assume that the equilibria for fixed $\mu$ exist for
$r\Gamma(\mu)<r_0$, we can choose the orientation of $v_0$ such that
$a_0>0$ and $a_2>0$. Thus, the reduced autonomous system \eqref{ydot}
has the form (for small $z\in\R$)
\begin{align}
\label{z dot1}
\dot{z} &= a_0(r\Gamma(\mu)-r_0) + a_2z^2+O(z^3), \\
\dot{\mu} &= \epsilon\Gamma(\mu)\mbox{.}
\label{mu dot}
\end{align} 
We zoom into the neighborhood of the maximal rate of change of
$\Gamma(\mu)$ (at $\mu_\mathrm{crit}$), $z=0$ and $r=r_0$ by
introducing rescaled variables and time
\begin{align*}
  \mu_\mathrm{old}&=\mu_\mathrm{crit}+\sqrt{\epsilon c_m}\mu_\mathrm{new}
  \mbox{,}&\mbox{where\quad}
   c_m&=\left[g_2r_0a_0a_2\right]^{-1/2}\mbox{,}\\
 z_\mathrm{old} &=\sqrt{\epsilon c_z}z_\mathrm{new}\mbox{,} &\mbox{where\quad}
 c_z&=\left[g_2r_0a_0/a_2^3\right]^{1/2}\mbox{,}\\
 r_\mathrm{old}&=r_0+\epsilon c_r r_\mathrm{new}\mbox{,}&\mbox{where\quad}
 c_r&=[g_2r_0/(a_0a_2)]^{1/2}\mbox{,}\\
 t_\mathrm{old}&=\sqrt{c_t/\epsilon}\,t_\mathrm{new}
  \mbox{,}&\mbox{where\quad}
c_t&=\left[g_2r_0a_0a_2\right]^{-1/2}\mbox{,}
\end{align*}
and expanding $\Gamma(\mu)$ near its unique maximum in
$\mu_\mathrm{crit}$ (recall from \eqref{Gamma properties} that
$\Gamma(\mu_{\mathrm{crit}})=1$ and
$g_2$ is defined as $-\Gamma''(\mu_\mathrm{crit})/2$):
\begin{align*}
  \Gamma(\mu_\mathrm{old}) &= \Gamma(\mu_{\mathrm{crit}}+\sqrt{\epsilon c_m}\mu_\mathrm{new})= 1 -\epsilon g_2 c_m\mu_\mathrm{new}^2+o(\epsilon)\mbox{.}
\end{align*}
In these coordinates the extended system \eqref{z dot1}--\eqref{mu
  dot} can then be written as
\begin{align}
\label{z dot2}
\dot{z} &= z^2+r-\mu^2+o(1)\mbox{,}\\
\dot{\mu} &= 1+O(\epsilon)
\label{m dot}
\end{align}
Orbits that stay close to the family of equilibria of \eqref{lin
  system} uniformly for all $\epsilon\to0$ are perturbations of orbits
that exist for all time in the limiting system of \eqref{z
  dot2}--\eqref{m dot} for $\epsilon=0$. We have $3$ cases for
\eqref{z dot2}--\eqref{m dot} with $\epsilon=0$. In all $3$ cases
there exists a unique \emph{globally defined} orbit $z(\mu)$ that exists for all times.
\begin{itemize}
\item[$r<1$] (Tracking): the globally defined orbit $z(\mu)$ has the
  limiting behavior $z(\mu)+|\mu|\to0$ for $\mu\to\pm\infty$ and is
  stable forward in time. All orbits starting with $z<0$ and $\mu\ll-1$
  converge to the globally defined orbit.
\item[$r=1$] (Critical): the globally defined orbit is $z(\mu)=\mu$.
\item[$r>1$] (Escape):  the globally defined orbit $z(\mu)$ has the
  limiting behavior $z(\mu)-|\mu|\to0$ for $\mu\to\pm\infty$ and is
  stable backward in time. All orbits starting with $z<0$ and $\mu\ll-1$
  diverge to $+\infty$ in finite time after $\mu>-\sqrt(r)$.
\end{itemize}
In the original coordinates the rescaled parameter $r=1$ equals the
first order expansion for the critical rate $r_c(\epsilon)$
\begin{equation*}
r_c(\epsilon) = r_0 + \epsilon\sqrt{-\frac{r_0\Gamma''(\mu_{\mathrm{crit}})}{2a_0a_2}}.
\end{equation*}
If white noise of variance $\sigma^2$ is added to \eqref{z dot1}:
\begin{displaymath}
  \d {z} = [a_0(r\Gamma(\mu)-r_0) + a_2z^2+O(z^3)]\d t+\sigma \d W_t\mbox{,}
\end{displaymath}
then $\sigma^2$ needs to be of the scale $\epsilon^{3/2}$. If
$\sigma^2=2D[\epsilon c_n]^{3/2}$ with
$c_n=(g_2r_0a_0)^{1/2}a_2^{-5/6}$ then the rescaled equation for $z$ is
\begin{displaymath}
  \d z = [z^2+r-\mu^2+o(1)]\d t+\sqrt{2D}\d W_t\mbox{.}
\end{displaymath}

\paragraph*{Limitation of single-mode approximation}
We use the saddle-node normal form \eqref{z dot2}--\eqref{m dot} to provide insight into the limitation of the single-mode approximation in reference to the maximal ramp speed $r$. The unique globally defined orbit $\tilde{z}(t)$ of \eqref{z dot2}--\eqref{m dot} for $\epsilon = 0$ is used to change to a co-moving coordinate system
\begin{equation*}
y(t)=z(t)-\tilde{z}(t),
\end{equation*}
with respect to $y$ such that \eqref{z dot2}--\eqref{m dot} has the form
\begin{align}
\nonumber
\dot{z} &= z^2 + r - t^2, \\
\nonumber
\dot{\tilde{z}} + \dot{y} &= (\tilde{z} + y)^2 + r - t^2, \\
\nonumber
\dot{y} &= \tilde{z}^2 + 2\tilde{z}y + y^2 + r - t^2 - \dot{\tilde{z}}, \\
\dot{y} &= y^2 + 2\tilde{z}y = y(y+2\tilde{z}).
\label{Comoving ODE}
\end{align} 
The potential $U(y,t)=-\int y(y+2\tilde{z})\,\mathrm{d}y$ in the new coordinate system has a well at $y=0$ and a hill top at $y=-2\tilde{z}$ for $\tilde{z}<0$, but a hill top at $y=0$ (and a well at $y=-2\tilde{z}$) for $\tilde{z}>0$. 

If white noise is added to \eqref{Comoving ODE}, the error of the single-mode approximation is small provided $|\dot{\tilde{z}}(t)|$ is small, which is the case for $t$ close to the maximum of $\tilde{z}$. If $\tilde{z}(t)<0$ for all $t$ then escape is most likely to occur close to time $t_{\max}$, the maximum of $\tilde{z}(t)$ since the potential barrier is at its lowest. The mode approximation error at times close to $t_{\max}$ is small because $|\dot{\tilde{z}}(t)|$ is small, as discussed in the paper.     

However, if $\tilde{z}(t)>0$ for a range of $t$ then escape occurs with a non-small probability at times when $|\dot{\tilde{z}}(t)|$ is not small. We identify that a maximal ramp speed $r\approx 0.59$ corresponds to $\max\{\tilde{z}(t):t\in[t_0,T_{\mathrm{end}}]\} = 0$. Therefore, in the original coordinates the single-mode approximation fails for some maximal ramp speeds $r$ less than $r_c$, and in particular, we consider only for parameters up to

\begin{equation*}
r = r_0 + \frac{1}{2}\epsilon c_r,
\end{equation*} 

\noindent where in the paper the constant $r_1$ is the same as $c_r$.  

%\noindent Note that the same critical rate is obtained for both orders of $m$. 

\section{Phase planes of shifted system}
\label{app:phase planes}

In this section we consider all qualitatively different phase planes for the shifted slow-fast autonomous system
\begin{align}
\label{SN ydot}
\dot{y} &= f(y) + br\Gamma(\mu), \\
\dot{\mu} &= \epsilon\Gamma(\mu).
\label{SN mudot}
\end{align}
%\noindent where $\mu$ is considered as a ramped shift, which has a maximal speed $r$. The direction of the shift is determined by $b$ and $\Gamma(\mu)$ satisfies the following properties
%
%\begin{equation*}
%  \begin{split}
%    \Gamma(0)&=\Gamma(1)=0\mbox{,\
%    }\max\{\Gamma(\mu):\mu\in[0,1]\}=1\mbox{,\ and\ }\\
%    \Gamma(\mu)&>0\mbox{\ for all $\mu\in(0,1)$.}
%  \end{split}
%\end{equation*}

\begin{figure}[h!]
        \centering
        \subcaptionbox{}[0.45\linewidth]
                {\includegraphics[scale = 0.3]{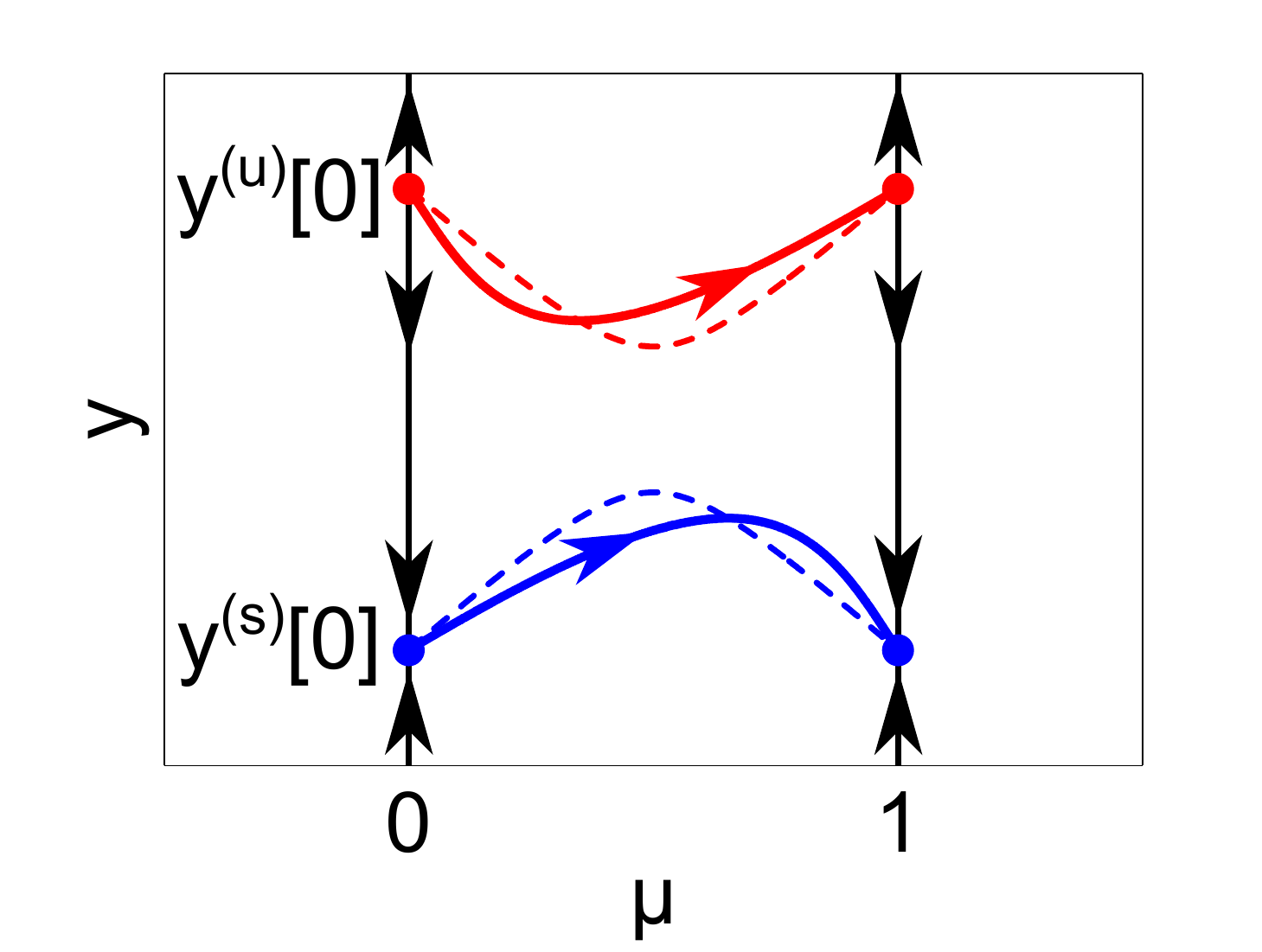}}
        \hfill
        \subcaptionbox{}[0.45\linewidth]
                {\includegraphics[scale = 0.3]{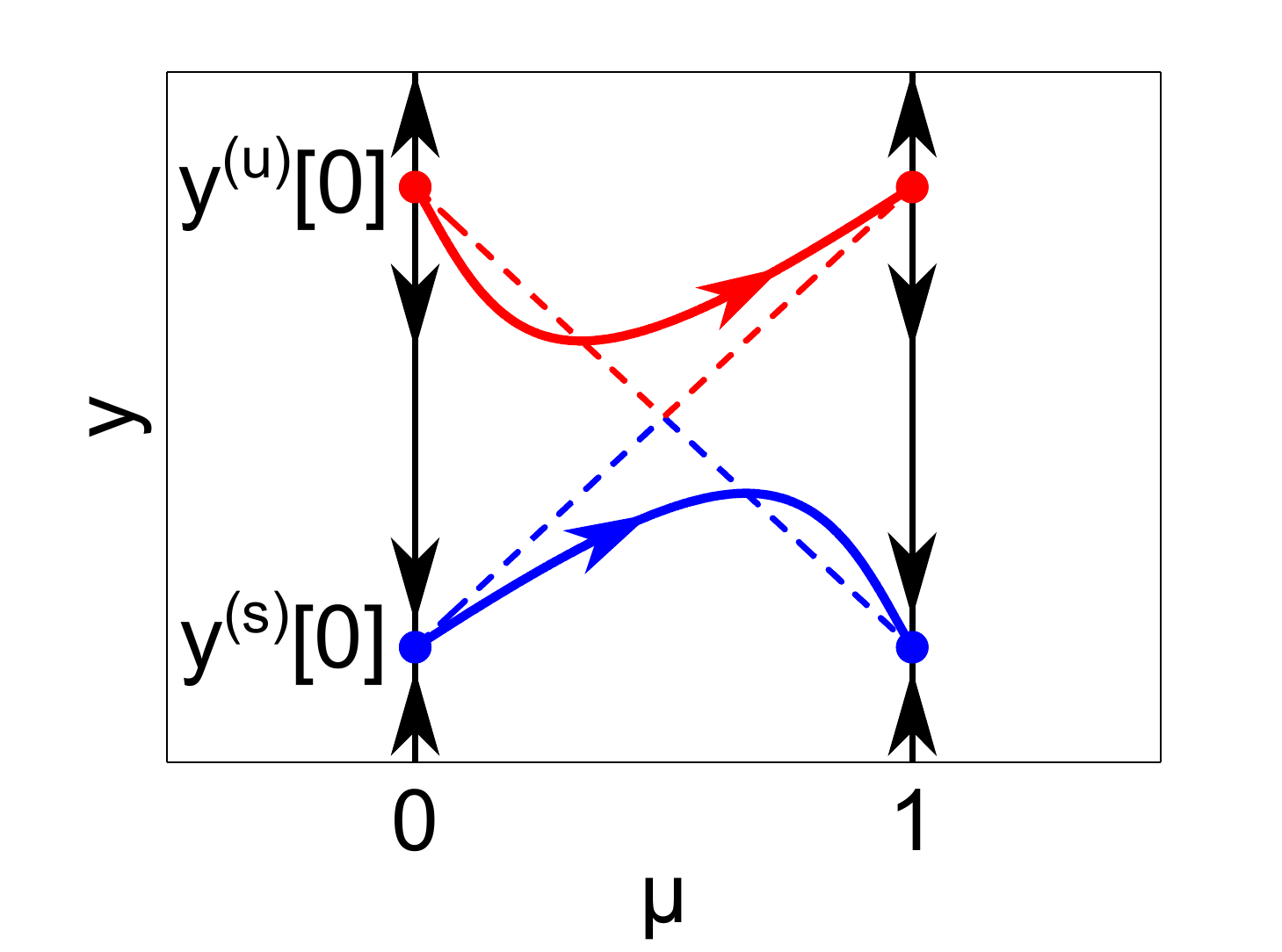}}
        ~ 
        \\
        \subcaptionbox{}[0.45\linewidth]
                {\includegraphics[scale = 0.3]{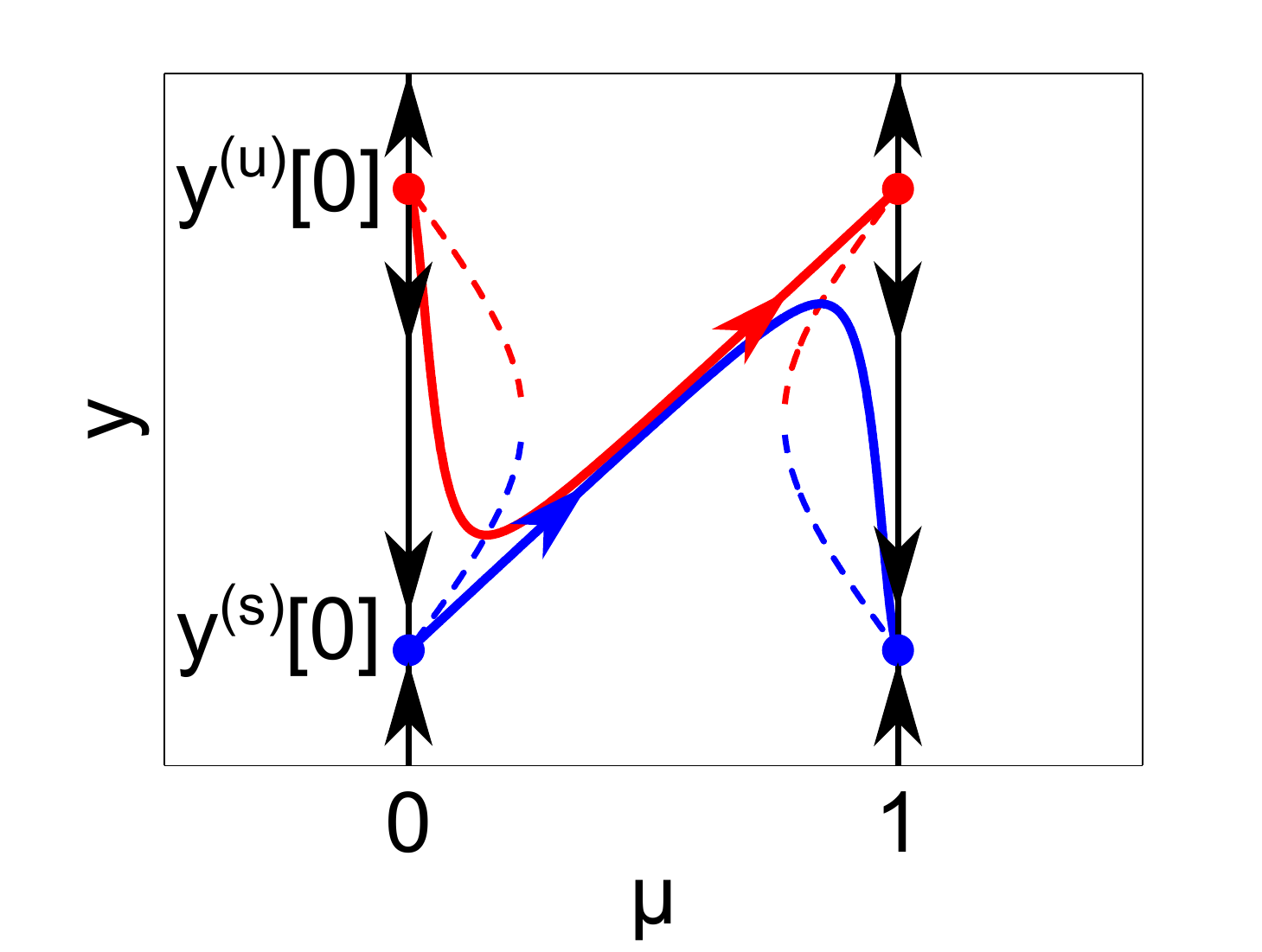}}
        \hfill
        \subcaptionbox{}[0.45\linewidth]
                {\includegraphics[scale = 0.3]{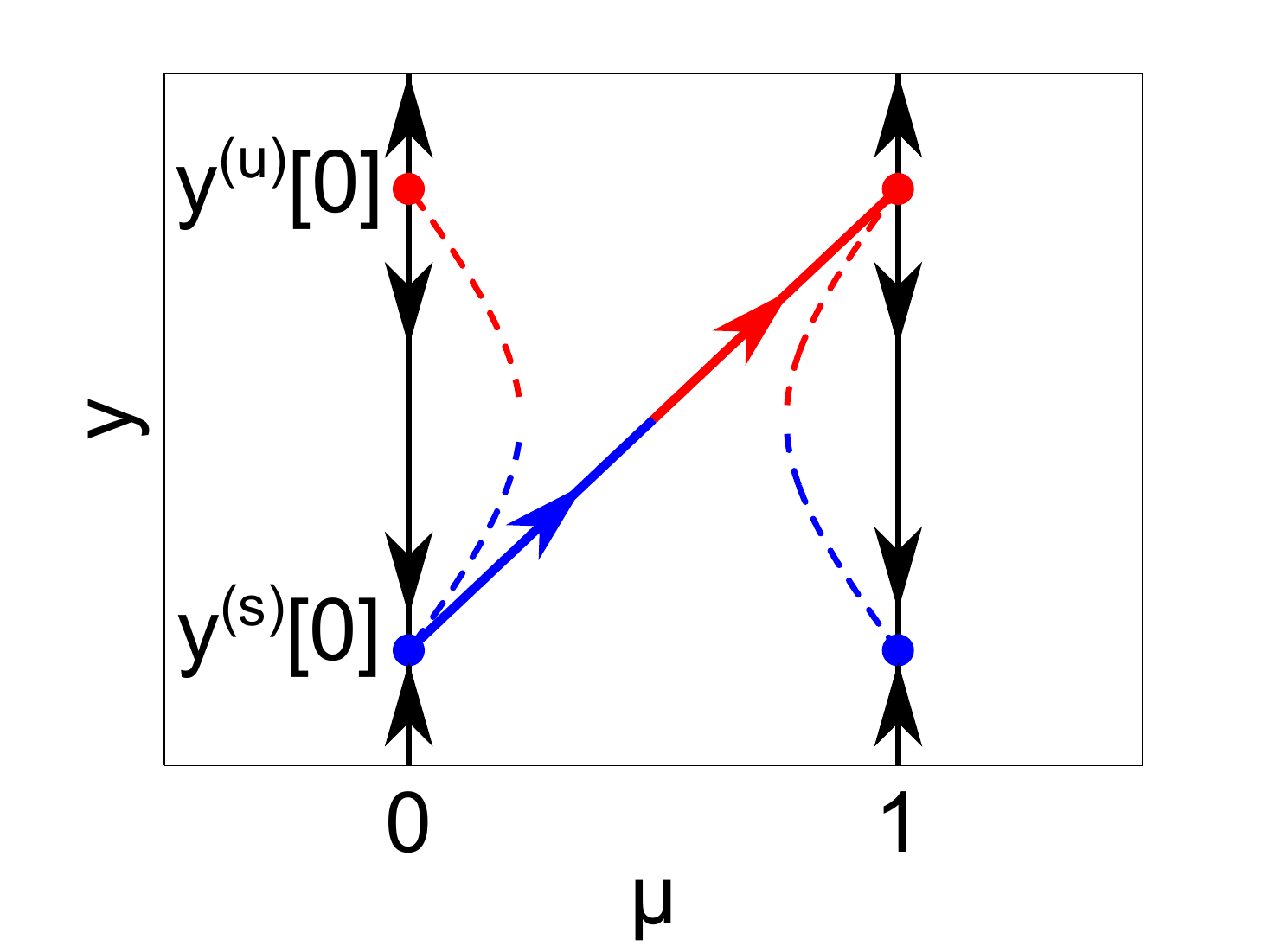}}
        ~ 
        \\
        \subcaptionbox{}[0.45\linewidth]
                {\includegraphics[scale = 0.3]{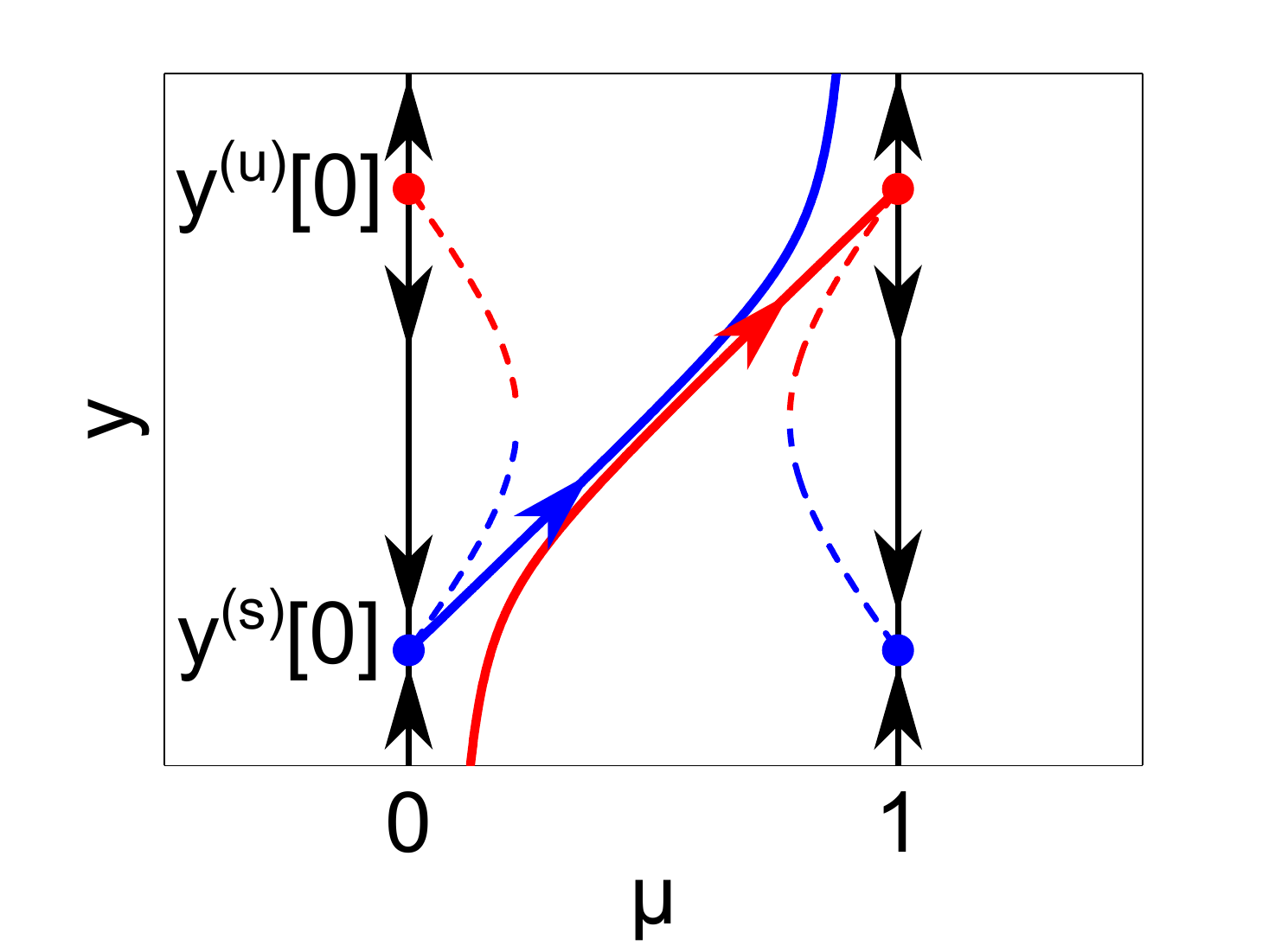}}
        \hfill
        \subcaptionbox{}[0.45\linewidth]
                {\includegraphics[scale = 0.3]{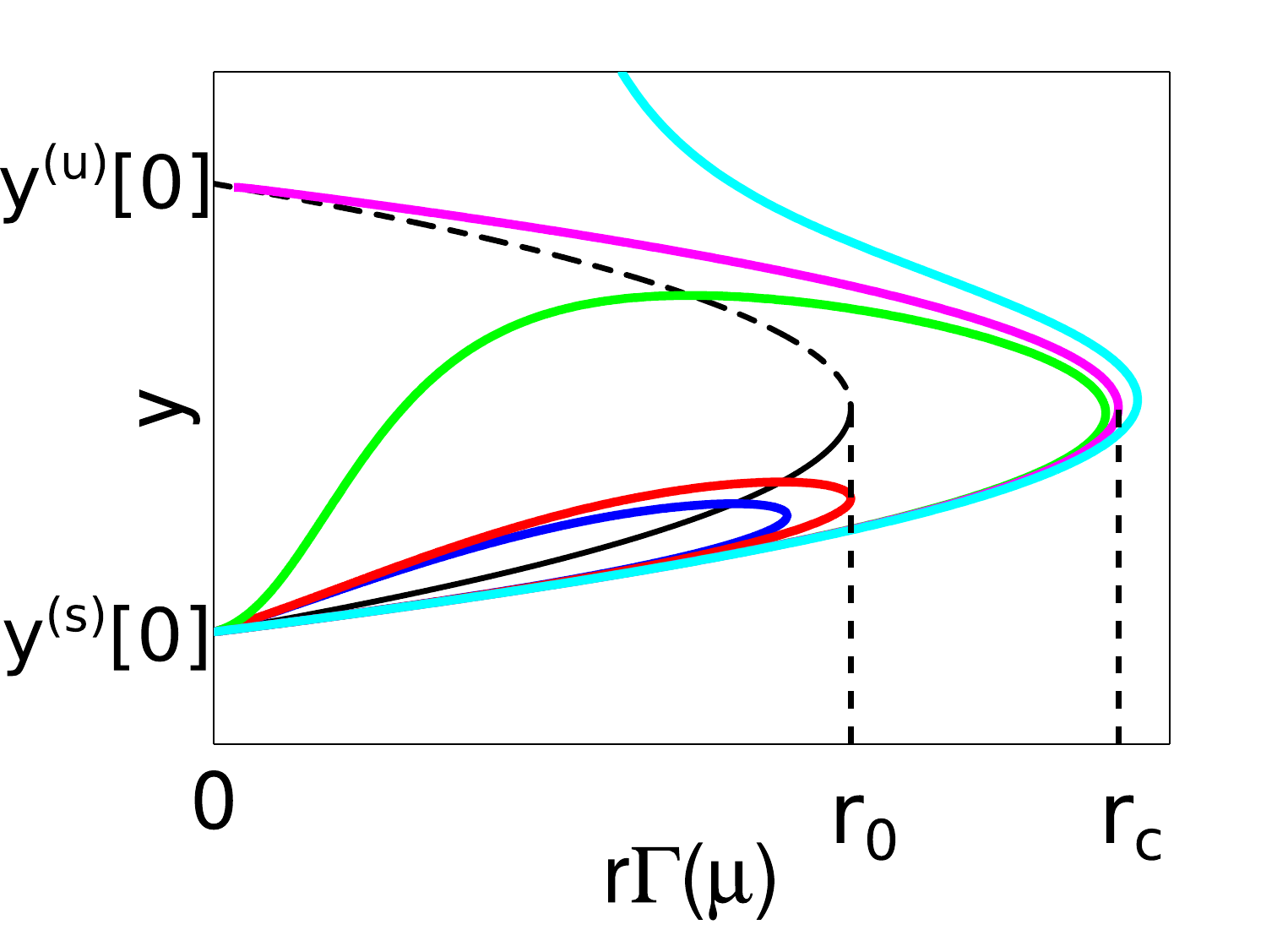}}
        ~ 
        \caption{(a)-(e) Qualitatively different phase planes of system \eqref{SN ydot}--\eqref{SN mudot} dependent on the value of the maximal ramp speed $r$ in relation to the saddle-node bifurcation at $r\Gamma(\mu)=r_0$ and the critical rate $r_c$ that induces tipping (see panel (f)). (a) $r<r_0<r_c$, (b) $r=r_0<r_c$, (c) $r_0<r<r_c$, (d) $r_0<r=r_c$, (e) $r_0<r_c<r$. Blue (lower) and red (upper) dashed lines are the stable and unstable branches of equilibria in the limit $\epsilon = 0$. Solid blue and red curves represent the unstable manifold of the saddle $(y,\mu)=(y^\s[0],0)$ and the stable manifold of the saddle $(y,\mu)=(y^\u[0],1)$ respectively. %Panels (a)-(c) correspond to tracking where the system \eqref{SN ydot}--\eqref{SN mudot} does not tip ($r<r_c$). Panel (d) displays the saddle-to-saddle connection for when $r=r_c$. Panel (e) is for $r>r_c$ which causes system \eqref{SN ydot}--\eqref{SN mudot} to tip. 
Panel (f) provides the bifurcation diagram of \eqref{SN ydot}--\eqref{SN mudot} in the $(r\Gamma(\mu),y)$ - plane. Black solid curve denotes branch of stable equilibria and black dashed curve branch of unstable equilibria. Superimposed on top are trajectories of a realisation starting arbitrarily close to $(y,\mu)=(y^\s[0],0)$ for panel (a) given in blue, (b) red, (c) green, (d) pink, (e) light blue (listed in increasing values of $r$). Parameter: $b=1$.}\label{app:Phase planes}
\end{figure}

\noindent Again we assume that in the limit $\epsilon = 0$, \eqref{SN ydot} has a saddle-node bifurcation at $r\Gamma(\mu)=r_0>0$, $y=y_0$ with a stable branch $y^{(s)}[r\Gamma(\mu)]$ and an unstable branch $y^{(u)}[r\Gamma(\mu)]$ of equilibria emerging for $r\Gamma(\mu)\in[0,r_0]$. The properties \eqref{Gamma properties} of $\Gamma$ imply that there exists three distinct sets of equilibria branches depending on the value of $r$ in relation to $r_0$. These are presented in the $(\mu,y)$ - phase plane in panel (a) for $r<r_0$, panel (b) for $r=r_0$ and panels (c)-(e) for $r>r_0$ in Figure \ref{app:Phase planes}. We will discuss the differences between panels (c)-(e) but initially we will just focus on the branches of stable (blue dashed) and unstable (red dashed) equilibria for $\epsilon=0$.   

For $r<r_0$ in the limit $\epsilon = 0$ there exists one stable equilibrium and one unstable equilibrium for all $\mu\in[0,1]$, see panel (a). This means there is a continuous branch of stable equilibria connecting $(y,\mu)=(y^\s[0],0)$ and $(y^\s[0],1)$ and likewise an unstable branch connecting $(y,\mu)=(y^\u[0],0)$ and $(y^\u[0],1)$. The phase portrait in panel (a) corresponds to the case $r<0$ for the rescaled $r$ in the rescaled system \eqref{z dot2}--\eqref{m dot}.

For $r=r_0$ the stable and unstable equilibria meet at $\mu_{\mathrm{crit}}\in(0,1)$ where $\Gamma(\mu_{\mathrm{crit}})=1$; see panel (b). In contrast to panel (a) the branches of equilibria approach at a linear rate and then also move away at a linear rate. This therefore means that the branches are non-differentiable at $\mu=\mu_{\mathrm{crit}}$ but the continuous connections still exist. 

In panels (c)-(e) the maximal ramp speed $r$ is greater than $r_0$ meaning that for $\mu$ in some interval $(\mu_l,\mu_u)\subset(0,1)$ centered around $\mu_{\mathrm{crit}}$ no equilibria exist in the limit $\epsilon = 0$. Instead two saddle-node bifurcations form at $\mu_l$ and $\mu_u$ with the branches of stable and unstable equilibria emerging for $\mu\in[0,\mu_l]$ and $\mu\in[\mu_u,1]$ to $y^\s[0]$ and $y^\u[0]$ respectively at $\mu=0$ and $1$. 

Let us now discuss the solid blue and red curves in panels (a)-(e) which denote the unstable manifold of $(y,\mu)=(y^\s[0],0)$ and the stable manifold of $(y,\mu)=(y^\u[0],1)$ respectively. Panels (a)-(c) depict the tracking scenario $(r<r_c)$ such that a connecting orbit exists between $y^\s[0]$ at $\mu=0$ and $y^\s[0]$ at $\mu=1$, the unstable manifold of  $(y,\mu)=(y^\s[0],0)$. The stable manifold of $(y,\mu)=(y^\u[0],1)$ acts as a separatrix where all solutions below the manifold converge to the stable node $(y,\mu)=(y^\s[0],1)$ and all those above escape to infinity. For panels (a) and (b), solutions $(y(t),\mu(t))$ starting close to $(y^\s[0],1)$ stay close to the stable equilibrium branch $y^\s[r\Gamma(\mu(t))]$ for all $t$ (distance goes to $0$ as $\epsilon\to 0$). Notice in panel (c) that although $r>r_0$ such that the saddle-node bifurcation is crossed for a small period of time the system \eqref{SN ydot}--\eqref{SN mudot} still does not tip because $r<r_c = r_0+\mathcal{O}(\epsilon)$ the critical rate. The scenario in panel (c) corresponds to the case $0<r<1$ for the rescaled $r$ in the rescaled system \eqref{z dot2}--\eqref{m dot}.

In panel (d) $r$ equals $r_c$, which creates a saddle-to-saddle connection from $(y^\s[0],0)$ to $(y^\u[0],1)$. Finally, in panel (e) $r>r_c$ which induces system \eqref{SN ydot}--\eqref{SN mudot} to tip. This means that initial conditions starting arbitrarily close to $(y^\s[0],0)$ go on to escape following the unstable manifold of $(y^\u[0],1)$.   Panels (d) and (e) correspond to the critical case $r=1$ (d) and the escape case $r>1$ (e) for the rescaled $r$ in the rescaled system \eqref{z dot2}--\eqref{m dot}.

We present in panel (f) the bifurcation diagram of \eqref{SN ydot}--\eqref{SN mudot} in the $(r\Gamma(\mu),y)$ - plane. The black solid and dashed curves give the stable and unstable branches of equilibria respectively. Superimposed on top are colored curves representing the trajectory for starting close to $(y,\mu)=(y^\s[0],0)$ for each of the scenarios in panels (a)-(e). The structure of the ramp means that solutions starting at $(r\Gamma(\mu),y)=(0,y^\s[0])$ approach the saddle-node at $r\Gamma(\mu)=r_0$ at a slow speed  initially which gets faster until $\mu=\mu_{\mathrm{crit}}$ is reached where $r\Gamma(\mu_{\mathrm{crit}}) = r$ (the turning points of the colored curves; recall, $\mu_{\mathrm{crit}}$ is where $\max\{\Gamma(\mu):\mu\in[0,1]\}=1$). The solutions then return back to $r\Gamma(1)=0$.

A prominent feature in all cases is that initially the trajectory appears to lag behind the stable equilibrium branch $y^\s[r\Gamma(\mu)]$ as $\Gamma(\mu(t))$ increases. The dark blue trajectory corresponds to panel (a) where $r<r_0$ and so in this scenario the saddle-node is not reached, which means the stable quasi-steady state is always present. Not evident from the bifurcation diagram is that when $\Gamma(\mu(t))$ reaches its maximum value the shift is at its fastest. Hence, the trajectory crosses the branch $y^\s[r\Gamma(\mu(t))]$ since the quasi-steady state changes direction quickly in comparison to the trajectory. The trajectory continues to lag until the ramp comes to a rest and the trajectory returns to $y^\s[0]$ at $r\Gamma(1)=0$. A similar pattern is observed for the red trajectory where $r=r_0$ (the scenario from panel (b)). The lag becomes more pronounced the closer $r\Gamma(\mu(t))$ gets to $r_0$ the saddle-node. Thus, only touching the saddle-node and not crossing it, combined with the speed of the shift, the trajectory will not escape and instead again cross the stable branch $y^\s[r\Gamma(\mu(t))]$ and converge back to $y^\s[0]$ at $r\Gamma(1)=0$.

The green curve shows the trajectory for $r_0<r<r_c$, which corresponds to panel (c). For $r\Gamma(\mu(t))>r_0$ the trajectory gives the impression that it begins to escape but only slowly. Though once again this is deceptive since the speed of the ramp is at its fastest during this phase and so it is only a short period of time before $\Gamma(\mu(t))=r_0$ for a second time. The quasi-steady states then begin to emerge but the trajectory is now above the unstable quasi-steady state. The trajectory continues to escape until it crosses the branch $y^\u[r\Gamma(\mu(t))]$. This implies that the quasi-steady state is shifting faster than the trajectory is escaping. Once across the unstable quasi-steady state the trajectory gets attracted back to $y^\s[0]$ at $r\Gamma(1)=0$ to complete the connecting orbit.

%\enlargethispage{\baselineskip}

Whereas, the pink trajectory ($r=r_c$, panel (d)) is shifted a little bit further such that the trajectory meets the $y^\u[r\Gamma(\mu(t))]$ branch only at the end of the ramp. The last scenario, $r>r_c$, the light blue trajectory is shifted sufficiently past the saddle-node bifurcation, such that the trajectory has enough time to escape before the system recovers.
\end{appendix}

\bibliography{Rate-induced_paper2}

\end{document}